%% file: NHIMs.tex
\documentclass[12pt]{iopart}
\usepackage{iopams}  
\usepackage{setstack}
\usepackage{amsfonts}

\usepackage{amsthm}
\usepackage{amssymb}
\usepackage{graphicx}
\usepackage{epsfig}
\usepackage{psfrag}

\newtheorem{theorem}{Theorem}
\newtheorem{assumption}[theorem]{Assumption}

\newtheorem{definition}[theorem]{Definition}
\newtheorem{lemma}[theorem]{Lemma}

\theoremstyle{remark}
\newtheorem{remark}[theorem]{Remark}
\newtheorem{example}[theorem]{Example}

\def\be{\begin{equation}}
\def\ee{\end{equation}}
\def\bea{\begin{eqnarray}}
\def\eea{\end{eqnarray}}

\def\Qset{\mathbb{Q}}
\def\Ocal{\mathcal{O}}

\def\eps{\varepsilon}

\def\hsm1{\hspace{-1mm}}

\begin{document}

\title[	Computer Assisted Proof for Normally Hyperbolic Invariant Manifolds]{Computer Assisted Proof for Normally Hyperbolic Invariant Manifolds}

\author{Maciej J. Capi\'{n}ski}

\address{Faculty of Applied Mathematics, AGH University of Science and Technology, al. Mickiewicza 30, 30-059 Krak\'{o}w, Poland}
\ead{mcapinsk@agh.edu.pl}

\author{Carles Sim\'o}

\address{Departament de Matem\`atica Aplicada i An\`alisi, Universitat de Barcelona, Gran Via 585, 08007 Barcelona, Spain}
\ead{carles@maia.ub.es}

\begin{abstract}
We present a topological proof of the existence of a normally hyperbolic
invariant manifold for maps. In our approach we do not require that the map is a
perturbation of some other map for which we already have an invariant manifold.
But a non-rigorous, good enough, guess is necessary. The required assumptions
are formulated in a way which allows for rigorous computer assisted
verification. We apply our method for a driven logistic map, for which
non-rigorous numerical simulation in plain double precision suggests the
existence of a chaotic attractor. We prove that this numerical evidence is false
and that the attractor is a normally hyperbolic invariant curve.
\end{abstract}

\maketitle

\input{01-intro.tex}
\input{02-setup.tex}

\input{03-geometric.tex}

\input{04-verif.tex}

\input{05-examples}
\input{06-bib.tex}

\end{document}

%% file: 01-intro.tex
\section{Introduction}

In this paper we give a proof of existence of normally hyperbolic invariant manifolds for maps. The construction is performed in the state space of the map. Assumptions needed for the proof are of twofold nature. First we require topological conditions which follow from suitable alignment of the coordinates (these are the so called covering relations). Next we require that our map satisfies cone conditions. The aim of the paper though is not to produce yet another proof of the normally hyperbolic invariant manifold theorem. Our aim is to produce a tool that can be applied in rigorous-computer-assisted proofs. To show the strength of our approach we apply our theorem to a driven logistic map introduced in \cite{BSV}. The considered map is such that standard numerical simulation gives evidence of a chaotic attractor. The example is a demonstration of the fact that one has to be careful with the arithmetics in simulations, since the numerical evidence of an attractor is false. The map in fact possesses a normally hyperbolic invariant curve. This is apparent when simulations are performed using multiple precision computations. The strength of our method lies in the fact that even for such an example, which defeats standard numerical simulations, we are able to produce a rigorous proof of existence of a normally hyperbolic invariant curve.

The approach to normally hyperbolic manifolds presented here is in the spirit of \cite{Ca} and \cite{CZ2}. In \cite{Ca} a topological proof of existence of invariant sets with normally hyperbolic type properties is given. In \cite{CZ2}  the result is extended to prove normally hyperbolic invariant manifolds. In both cases the proofs relied on assumptions that the first iterate of the map is well aligned with the stable and unstable manifolds. Similar approach was also used in \cite{CR} to give a proof of existence of a center manifold. The result in \cite{CR} is for ODEs and relies also on the fact that hyperbolic dynamics is uniform. The main difference between our paper and results mentioned above is that we assume that hyperbolic expansion and contraction aligns with the tangent spaces of the invariant manifolds after a suitable (possibly large) number of iterates of the map. This setting is more general, and also more typical for normal hyperbolicity.  

The paper is organized as follows. Section \ref{sec:setup} introduces basic notations used throughout the paper and provides a setup and an outline of our problem.  Section \ref{sec:geometric} contains a geometric construction of a normally hyperbolic manifold. We first give a construction of a "center-stable" manifold (the term "center-stable" refers to the normally hyperbolic invariant manifold union its associated
stable manifold; analogous terminology is used by us for the "center-unstable"
manifold). A center-unstable manifold is obtained using a mirror construction to the center-unstable manifold, by considering the inverse map. The intersection of the center-stable and center-unstable manifolds gives us the normally hyperbolic invariant manifold. In Section \ref{sec:ver-cond} we show how to verify assumptions of our theorems using local bounds on derivatives of the map. In Section \ref{sec:examples} we present our example of the driven logistic map and apply our method to it.          

%% file: 02-setup.tex


\section{Setup}\label{sec:setup}

We start by writing out some basic notations which we shall use throughout the paper. A notation $B_{i}(q,r)$ will stand for a ball of radius $r$ centered at $q$
in $\mathbb{R}^{i}.$ We will also use a notation $B_{i}=B_{i}(0,1)$. For a set
$A$ we will denote by $\overline{A}$ its closure, by $\mathrm{int\,}A$ its
interior and by $\partial A$ its boundary. For a function $f$ we will use a notation $\mathrm{dom}%
(f)$ to denote its domain. For points $p=(x,y)$ we shall use notation
$\pi_{x}(p),$ $\pi_{y}(p)$ to denote the projection onto the $x$ and $y$
coordinates respectively.

We now introduce the setup of our problem. Let $D$ and $\mathcal{U}$ be open subsets in $\mathbb{R}^{n}$ such that
$D\subset\mathcal{U}$. Let%

\[
f:\mathcal{U}\rightarrow\mathcal{U},
\]
be a diffeomorphism. Let $u,s,c\in\mathbb{N}$ be such that $u+s+c=n.$ We
assume that there exist a diffeomorphism%
\[
\phi:\mathcal{U}\rightarrow\phi(\mathcal{U})\subset\mathbb{R}^{u}%
\times\mathbb{R}^{s}\times\Lambda
\]
such that $\phi(\mathrm{cl\,}D)=D_{\phi}:=\overline{B}_{u}\times\overline{B}%
_{s}\times\Lambda$, and $\Lambda$ is a compact $c$ dimensional manifold
without boundary. We define $f_{\phi}:D_{\phi}\rightarrow\mathbb{R}^{u}%
\times\mathbb{R}^{s}\times\Lambda$ as
\[
f_{\phi}=\phi\circ f\circ\phi^{-1}.
\]
We assume that there exists a finite covering $\{U_{i}\}_{i\in I}$ of
$\Lambda$ and an atlas%
\[
\eta_{i}:\overline{U}_{i}\rightarrow\overline{B}_{c}.
\]
Throughout the work we will use a notation%
\[
\mathbf{B}=\overline{B}_{u}\times\overline{B}_{s}\times\overline{B}_{c}.
\]
For $i,j\in I$ we consider local maps $f_{ji}:\mathbf{B}\supset\mathrm{dom}%
(f_{ij})\rightarrow\mathbb{R}^{u}\times\mathbb{R}^{s}\times\overline{B}_{c}$
defined as%
\begin{eqnarray*}
f_{ij}  &  := &\tilde{\eta}_{j}\circ f_{\phi}\circ\tilde{\eta}_{i}^{-1}, \\
\tilde{\eta}_{i}  &  := &(\mathrm{id},\mathrm{id},\eta_{i})\qquad\text{for }i\in
I. 
\end{eqnarray*}
Note that the domain of $f_{ij}$ can be empty, and will usually be smaller
than $\mathbf{B}.$ The following graph depicts the above defined functions and
their mutual relations.
\[%
\begin{array}
[c]{ccc}%
D & \overset{f}{\rightarrow} & \mathcal{U}\\
\downarrow\phi &  & \downarrow\phi\\
\overline{B}_{u}\times\overline{B}_{s}\times\Lambda & \overset{f_{\phi}%
}{\rightarrow} & \mathbb{R}^{u}\times\mathbb{R}^{s}\times\Lambda\\
\quad\downarrow\tilde{\eta}_{i} &  & \quad\downarrow\tilde{\eta}_{j}\\
\mathbf{B} & \overset{f_{ji}}{\rightarrow} & \mathbb{R}^{u}\times
\mathbb{R}^{s}\times\overline{B}_{c}%
\end{array}
\]

Our task in this paper will be to find a normally hyperbolic invariant
manifold, together with its stable and unstable manifolds within the set $D$.

We will use the following notations for our coordinates: $x\in\mathbb{R}^{u},$
$y\in\mathbb{R}^{s},$ $\theta\in\overline{B}_{c},$ $\lambda\in\Lambda.$ The
coordinate $x$ will play the role of a globally unstable direction, and the
coordinate $y$ will play the role of a stable direction for the map $f_{\phi}$
(hence the superscripts $u$ and $s$, which stand for "unstable" and "stable"
respectively). The coordinate $\lambda$ will play the role of the central
direction, in which the global dynamics is weaker than in the stable and unstable
coordinates. The notation $\theta$ will also be used for the central
direction, but it will be reserved to denote the central coordinate in the
local coordinates; i.e. $\theta=\eta_{i}(\lambda)$ for some $\lambda\in
\Lambda$ and $i\in I$.

%% file: 03-geometric.tex


\section{Geometric approach to invariant manifolds}\label{sec:geometric}

In this section we give the construction of a normally hyperbolic invariant
manifold. The construction is performed in the state space of our map. It is
based on the assumptions of covering relations and cone conditions. We first
give an introduction to these tools in Section \ref{sec:cov-cc}. In Section
\ref{sec:cover-cc-norm-hyp} we formulate our assumptions on the map in terms
of covering relations and cone conditions, which will imply the existence of a
normally hyperbolic manifold. In Section \ref{sec:main-res} we show how to
construct a center-stable manifold of our map. The construction of a
center-unstable manifold follows from a mirror argument. The intersection of
center-stable and center-unstable manifolds gives us a $C^{0}$ normally
hyperbolic invariant manifold. Let us write explicitly that for a normally
hyperbolic manifold which does not have an associated stable manifold, the
center-stable manifold will be the the normally hyperbolic manifold itself.
Analogous statement holds also for center-unstable manifolds.

\subsection{Covering relations and cones}

\label{sec:cov-cc} Covering relations are topological tools used for proofs of
nontrivial symbolic dynamics of dynamical systems. The method is based on the
Brouwer fixed point index, and the setting is such that it allows for rigorous
numerical verification. The method has been applied in computer assisted
proofs for the H\'{e}non map, R\"{o}ssler equations \cite{Z}, \cite{CZ2},
Lorenz equations \cite{GaZ}, Chua circuit \cite{G} or Kuramoto-Shivashinsky
ODE \cite{W}, amongst others. The method is based on singling out a number of
regions, called h-sets, which have hyperbolic type properties. Using these
properties one can find orbits of the system, which shadow the h-sets along
their trajectories. The method of covering relations relies on the system
having expanding and contracting coordinates. In this section we generalize
covering relations to include also a central direction. The setup is similar
to that of \cite{Ca}, \cite{CZ}, but has been simplified. Our proofs are now
simpler and based only on continuity arguments. They no longer require the use
of degree theory, with little loss of generality.

For any $p=(x,y,\theta)\in\mathbf{B}$ and $r_{u},r_{s},r_{c}>0$ we introduce a
notation
\[
N(p,r_{u},r_{s},r_{c}):=\overline{B}_{u}(x,r_{u})\times\overline{B}%
_{s}(y,r_{s})\times\overline{B}_{c}(\theta,r_{c}).
\]
We define%
\begin{eqnarray*}
N^{-}  &  = & N^{-}(p,r_{u},r_{s},r_{c}):=\partial\overline{B}_{u}(x,r_{u}%
)\times\overline{B}_{s}(y,r_{s})\times\overline{B}_{c}(\theta,r_{c})\\
N^{+}  &  = & N^{+}(p,r_{u},r_{s},r_{c}) \\
&  := &\overline{B}_{u}(x,r_{u})\times\left(  (\mathbb{R}^{s}\times
\mathbb{R}^{c})\setminus(B_{s}(y,r_{s})\times B_{c}(\theta,r_{c}))\right)  .
\end{eqnarray*}
We assume that all boxes $N$ which we are going to consider here are contained
in $\mathbf{B}$. We will refer to a box $N$ as a \emph{ch-set}
(center-hyperbolic set) centered at $p$.

In following arguments we shall often consider different ch-sets. To keep
better track of our notations and to make our arguments more transparent we
shall stick to a convention that for two ch-sets $N_{1},N_{2}$ centered
respectively at $p_{1}=(x_{1},y_{1},\theta_{1})$ and $p_{2}=(x_{2}%
,y_{2},\theta_{2})$ we shall write%
\[
N_{i}=N_{i}(p_{i},r_{u}^{i},r_{s}^{i},r_{c}^{i}):=\overline{B}_{u}^{i}%
(x_{i},r_{u}^{i})\times\overline{B}_{s}^{i}(y_{i},r_{s}^{i})\times\overline
{B}_{c}^{i}(\theta_{i},r_{c}^{i})\quad\text{for }i=1,2.
\]
\begin{figure}[ptb]
\begin{center}
\includegraphics[
width=3in
]{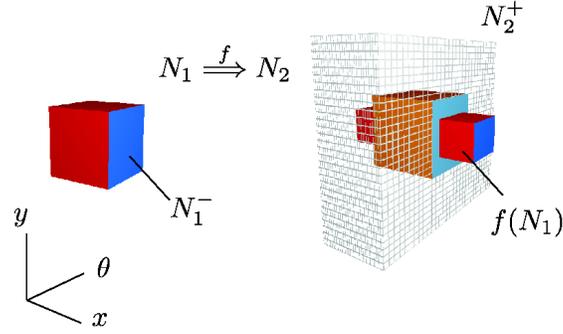}
\end{center}
\caption{A ch-set $N_{1}$ covering a ch-set $N_{2}.$}%
\label{fig:covering1}%
\end{figure}

\begin{definition}
\label{def:f-covers}Let $g:\mathbf{B}\rightarrow\mathbb{R}^{u}\times
\mathbb{R}^{s}\times\overline{B}_{c}$ be a continuous function. Let
$p_{i}=(x_{i},y_{i},\theta_{i})$ for $i=1,2$ and let $N_{1}$, $N_{2}$ be two
ch-sets in $\mathbf{B}$ centered at $p_{1}$ and $p_{2}$ respectively. We say
that $N_{1}$ $g$-\emph{covers} $N_{2}$ if%
\begin{eqnarray}
g(p_{1})   \in  \mathrm{int}(N_{2}),\label{eq:covering-cond-1}\\
\pi_{x}(g(N_{1}^{-}))\cap\overline{B}_{u}^{2}(x_{2},r_{u}^{2})  =  \emptyset,\label{eq:covering-cond-2}\\
g(N_{1})\cap N_{2}^{+}   =  \emptyset. \label{eq:covering-cond-3}%
\end{eqnarray}
In such case we shall write $N_{1}\overset{g}{\Longrightarrow}N_{2}.$
\end{definition}

\begin{remark}
Definition \ref{def:f-covers} is a simplified definition of a covering
relation. More general versions can be found in \cite{GaZ}, \cite{GiZ},
\cite{Z} in the setting of hyperbolicity, or in \cite{Ca}, \cite{CZ} in a
setting when additionally a central direction is included.
\end{remark}

\begin{figure}[ptb]
\begin{center}
\includegraphics[
height=2.6in
]{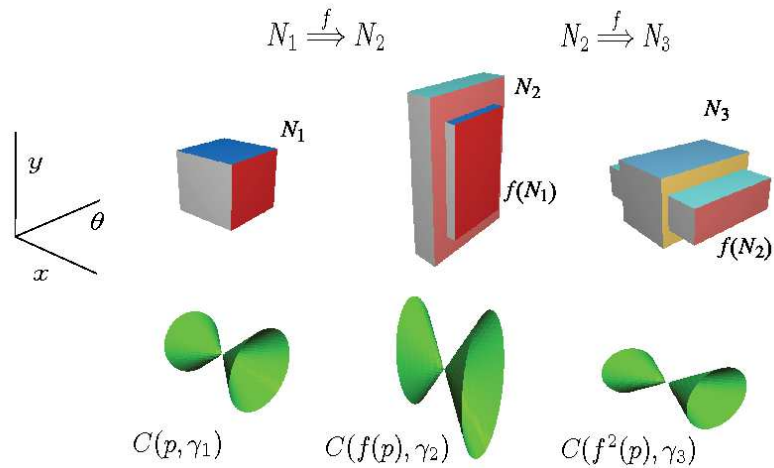}
\end{center}
\caption{Covering relations for two iterates of a map $f$. For the second
iterate of the map the coordinate $x$ is expanding and $y$ is contracting (for
the first iterate of $f$ they are not). The fact that expansion in $x$ is
stronger than expansion in $\theta$ is visible from the fact that the cones
$C(f^{2}(p),\gamma_{3})$ are "tighter" than cones $C(p,\gamma_{1})$.}%
\label{fig:covering-steps}%
\end{figure}

For $\gamma=(a,b,c)\in\mathbb{R}^{3},$ and $q=(x,y,\theta)\in\mathbb{R}%
^{u}\times\mathbb{R}^{s}\times\mathbb{R}^{c}$ we define%
\[
Q_{\gamma}:\mathbb{R}^{u}\times\mathbb{R}^{s}\times\mathbb{R}^{c}%
\rightarrow\mathbb{R}%
\]%
\begin{equation}
Q_{\gamma}(q):=a\left\Vert x\right\Vert ^{2}+b\left\Vert y\right\Vert
^{2}+c\left\Vert \theta\right\Vert ^{2}. \label{eq:Q-formula}%
\end{equation}
If $a>0$ $b,c<0,$ then for $p\in\mathbb{R}^{u}\times\mathbb{R}^{s}%
\times\mathbb{R}^{c}$ we will refer to%
\[
C(p,\gamma):=\{q:Q_{\gamma}(p-q)\geq0\}
\]
as a \emph{horizontal cone} centered at $p$ (see Figure
\ref{fig:covering-steps}).

\begin{definition}
Let $N$ be a ch-set and $\gamma=(a,b,c)$ be such that $a>0,$ $b,c<0.$ We will
refer to a pair $(N,\gamma)$ as a \emph{ch-set with cones}.
\end{definition}

\begin{definition}
\label{def:horizontal-disc}Let $(N,\gamma)=(N((x,y,\theta),r_{u},r_{s}%
,r_{c}),\gamma)$ be a ch-set with cones. A continuous function $\mathbf{h}%
:\overline{B}_{u}(x,r_{u})\rightarrow N$ is called a \emph{horizontal disc} in
$(N,\gamma),$ iff $\pi_{x}\mathbf{h}(x)=x$ and for any $x^{\ast},x^{\ast\ast
}\in\overline{B}_{u}(x,r_{u}),$%
\begin{equation}
Q_{\gamma}(\mathbf{h}(x^{\ast})-\mathbf{h}(x^{\ast\ast}))\geq0,
\label{eq:b-horizontal-disc-ineq}%
\end{equation}

\end{definition}

\begin{lemma}
\label{lem:hor-disc-image}Let $N_{i}=N_{i}((x_{i},y_{i},\theta_{i}),r_{u}%
^{i},r_{s}^{i},r_{c}^{i})$ for $i=1,2$ and let $(N_{1},\gamma_{1})$,
$(N_{2},\gamma_{2})$ be two ch-sets with cones. Assume that%
\begin{equation}
N_{1}\overset{g}{\Longrightarrow}N_{2} \label{eq:N1-N2-covering}%
\end{equation}
and that for any $q^{\ast},q^{\ast\ast}\in N_{1}$ such that $q^{\ast}\neq
q^{\ast\ast}$ and $Q_{\gamma_{1}}(q^{\ast}-q^{\ast\ast})\geq0$ we have%
\begin{equation}
Q_{\gamma_{2}}(g(q^{\ast})-g(q^{\ast\ast}))>0. \label{eq:cone-cond-N1-N2}%
\end{equation}
If $\mathbf{h}_{1}$ is a horizontal disc in $(N_{1},\gamma_{1})$ then there
exists a horizontal disc $\mathbf{h}_{2}$ in $(N_{2},\gamma_{2})$ such that
$g(\mathbf{h}_{1}(\overline{B}_{u}^{1}(x_{1},r_{u}^{1})))\cap N_{2}%
=\mathbf{h}_{2}(\overline{B}_{u}^{2}(x_{2},r_{u}^{2})).$
\end{lemma}

\begin{proof}
Without loss of generality we assume that $p_{1}=p_{2}=0$ and that $r_{\kappa
}^{i}=1$ for $i=1,2$ and $\kappa\in\{u,s,c\}$. In other words we assume that
for $i=1,2$%
\[
N_{i}=\overline{B}_{u}^{i}\times\overline{B}_{s}^{i}\times\overline{B}_{c}%
^{i}=\overline{B}_{u}(0,1)\times\overline{B}_{s}(0,1)\times\overline{B}%
_{c}(0,1).
\]

Let $\gamma_{i}=(a_{i},b_{i},c_{i})$ for $i=1,2$ and let $\mathbf{h}$ be any
horizontal disc in $N_{1}.$ Then by (\ref{eq:Q-formula}),
(\ref{eq:b-horizontal-disc-ineq}) and (\ref{eq:cone-cond-N1-N2}) for $x^{\ast
},x^{\ast\ast}\in\overline{B}_{u}^{1},$ $x^{\ast}\neq x^{\ast\ast}$%
\begin{equation}
a_{2}\left\Vert \pi_{x}g(\mathbf{h}(x^{\ast}))-\pi_{x}g(\mathbf{h}(x^{\ast
\ast}))\right\Vert ^{2}\geq Q_{\gamma_{2}}(g(\mathbf{h}(x^{\ast}%
))-g(\mathbf{h}(x^{\ast\ast})))>0,\label{eq:cone-cond-proof1}%
\end{equation}
which means that $\pi_{x}\circ g\circ\mathbf{h}$ is a monomorphism.

Using a notation $\mathbf{h}_{1}(x)=(x,h_{1}(x))\in\overline{B}_{u}^{i}%
\times(\overline{B}_{s}^{i}\times\overline{B}_{c}^{i}),$ for $\alpha\in
\lbrack0,1],$ we define a family of horizontal discs $\mathbf{h}_{\alpha
}(x)=(x,\alpha h_{1}(x)).$ Let $F_{\alpha}:\overline{B}_{u}^{1}\rightarrow
\mathbb{R}^{u}$ be a continuous family of functions defined as%
\[
F_{\alpha}(x):=\pi_{x}\circ g\circ\mathbf{h}_{\alpha}(x).
\]
We shall show that $\overline{B}_{u}^{2}\subset F_{1}(B_{u}^{1}).$ Functions
$F_{\alpha}$ are monomorphisms, hence sets $A_{\alpha}:=F_{\alpha}(B_{u}^{1})$
are homeomorphic to balls in $\mathbb{R}^{u}$; moreover $\partial A_{\alpha
}=F_{\alpha}(\partial B_{u}^{1})$. By Definition \ref{def:horizontal-disc} of
a horizontal disc, $\mathbf{h}_{\alpha}(\partial B_{u}^{1})\subset N_{1}^{-}$.
From assumption (\ref{eq:N1-N2-covering}), by conditions
(\ref{eq:covering-cond-1}), (\ref{eq:covering-cond-2})%
\begin{eqnarray}
\pi_{x}g(0)   \in \overline{B}_{u}^{2},\label{eq:g-zero-B2}\\
\partial A_{\alpha}\cap\overline{B}_{u}^{2}   \subset  F_{\alpha}(N_{1}
^{-})\cap\overline{B}_{u}^{2}=\emptyset.\label{eq:A-alpha-bd}%
\end{eqnarray}
From the fact that $0\in B_{u}^{1}$
\begin{equation}
F_{0}(0)\in F_{0}(B_{u}^{1})=A_{0}.\label{eq:inter-ne1}%
\end{equation}
Since $\mathbf{h}_{0}(0)=0,$ by (\ref{eq:g-zero-B2})
\begin{equation}
F_{0}(0)=\pi_{x}\circ g\circ\mathbf{h}_{0}(0)=\pi_{x}g(0)\in\overline{B}%
_{u}^{2},\label{eq:inter-ne2}%
\end{equation}
From (\ref{eq:inter-ne1}), (\ref{eq:inter-ne2}) follows that $A_{0}%
\cap\overline{B}_{u}^{2}\neq\emptyset$. This by (\ref{eq:A-alpha-bd}) implies
that $\overline{B}_{u}^{2}\subset A_{0}$. By continuity of $F_{\alpha}$ with
respect to $\alpha$ this means that $\overline{B}_{u}^{2}\subset A_{\alpha}$
for all $\alpha\in\lbrack0,1].$ In particular $\overline{B}_{u}^{2}\subset
A_{1}=F_{1}(B_{u}^{1}).$

Since $F_{1}$ is a monomorphism and $\overline{B}_{u}^{2}\subset F_{1}%
(B_{u}^{1}),$ for any $v\in\overline{B}_{u}^{2}$ there exists a unique
$x=x(v)\in B_{u}^{1}$ such that $F_{1}(x)=v.$ We define $\mathbf{h}%
_{2}(v)=(v,h_{2}(v)):=(v,\pi_{y,\theta}\circ g\circ\mathbf{h}_{1}(x(v))).$ For
any $v^{\ast}\neq v^{\ast\ast}$, $v^{\ast},v^{\ast\ast}\in\overline{B}_{u}%
^{2}$, by (\ref{eq:b-horizontal-disc-ineq}) and (\ref{eq:cone-cond-N1-N2}) we
have%
\begin{eqnarray*}
Q_{\gamma_{2}}\left(  \mathbf{h}_{2}(v^{\ast})-\mathbf{h}_{2}(v^{\ast
})\right)   &  = & Q_{\gamma_{2}}(g\circ\mathbf{h}_{1}(x(v^{\ast}))-g\circ
\mathbf{h}_{1}(x(v^{\ast\ast})))\\
&  > &Q_{\gamma_{1}}(\mathbf{h}_{1}(x(v^{\ast}))-\mathbf{h}_{1}(x(v^{\ast\ast
})))\\
&  > &0.
\end{eqnarray*}
Since $Q_{\gamma_{2}}\left(  \mathbf{h}_{2}(v^{\ast})-\mathbf{h}_{2}%
(v^{\ast\ast})\right)  >0$%
\begin{eqnarray*}
a_{2}\left\Vert v^{\ast}-v^{\ast\ast}\right\Vert  \\
 >  -b_{2}\left\Vert \pi_{y}\left(  h_{2}(v^{\ast})-h_{2}(v^{\ast\ast})\right)  \right\Vert ^{2}%
-c_{2}\left\Vert \pi_{\theta}\left(  h_{2}(v^{\ast})-h_{2}(v^{\ast\ast
})\right)  \right\Vert ^{2}\\
\geq  \min(-b_{2},-c_{2})\left\Vert h_{2}(v^{\ast})-h_{2}(v^{\ast\ast
})\right\Vert ^{2},
\end{eqnarray*}
and therefore $\mathbf{h}_{2}$ is continuous.
\end{proof}

\begin{remark}
Let us note that since we have freedom of choice of the radii $r_{u},r_{s}$
and $r_{c}$ it is not necessary for $x$ to be expanding, $y$ to be contracting
and $\theta$ to have weaker dynamics for each single iterate of the map. In
Figure \ref{fig:covering-steps} we have a sketch of a situation in which $x$
becomes expanding and $y$ contracting after a second iterate. In Figure
\ref{fig:covering-steps} the coordinate $\theta$ is expanding. It will turn
out that such a scenario is acceptable for us and can be dealt with by
increasing $r_{c}$ for successive iterates.
\end{remark}

\subsection{Covering relations and cone conditions for normal hyperbolicity}

\label{sec:cover-cc-norm-hyp}

In this section we formulate our assumptions which will imply the existence of
a normally hyperbolic manifold. The assumptions are in terms of covering
relations and cones and are in the spirit of \cite{CZ}. There are two major
differences though. The first is that assumptions used in \cite{CZ} required
the system to have uniform expansion and uniform contraction for the first
iterate of the map. Here we set up our coordinates in the directions of
\emph{global} contraction and \emph{global} expansion. In the setting of
normal hyperbolicity the coordinates of global contraction and expansion need
not be contracting and expanding for the first iterates of the map. What is
important is that they dominate after a sufficiently large numbers of
iterates, in other words, that the Lyapunov exponents are negative or
positive, respectively. We set up our assumptions so that they allow for such
setting. The second difference is that our setup has been significantly
simplified with comparison to \cite{CZ}. This resulted in a slight loss of
generality (we do not formulate our assumptions in terms of vector bundles as
in \cite{CZ}) but we need to consider fewer assumptions.

Let $1>R>\rho,r>0.$ Assume that there exists a finite sequence of points
$\boldsymbol{\lambda}_{k}\in\Lambda,$ $k\in\mathbb{N}$ such that for any $k$
the set $I(k)=\{i:B_{c}(\eta_{i}(\boldsymbol{\lambda}_{k}),\rho)\subset
B_{c}(0,R)\}$ is not empty. What is more, assume that there exists a set
$J\subset\{(i,k)|i\in I(k)\}$ such that $\Lambda\subset\bigcup_{(i,k)\in
J}\eta_{i}^{-1}(B_{c}(\eta_{i}(\boldsymbol{\lambda}_{k}),\rho)).$ For points
$(i,k)\in J$ we define sets
\[
M_{i,k}:=\overline{B}_{u}(0,r)\times\overline{B}_{s}(0,r)\times\overline
{B}_{c}(\eta_{i}(\boldsymbol{\lambda}_{k}),\rho).
\]

\begin{figure}[ptb]
\begin{center}
\includegraphics[
height=2in
]{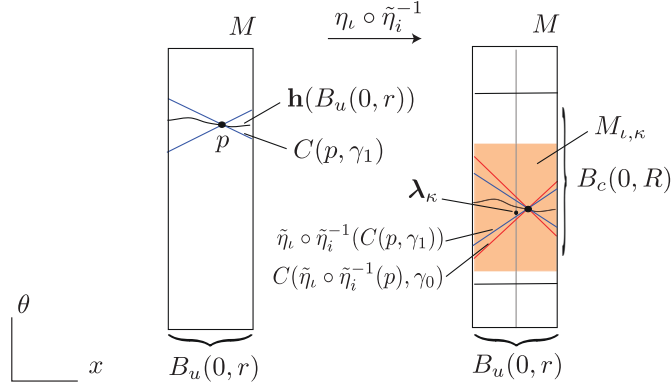}
\end{center}
\caption{The change of coordinates $\tilde{\eta}_{\iota}\circ\tilde{\eta}%
_{i}^{-1}$, a horizontal disc $\mathbf{h,}$ and the cones given by
$\boldsymbol{\gamma}_{0}$ and $\boldsymbol{\gamma}_{1}$ in different local
coordinates. Here, for simplicity, the stable coordinate is neglected $s=0$. }%
\label{fig:M-assumptions}%
\end{figure}

We will need to assume that the points $\boldsymbol{\lambda}_{k}$ are
sufficiently close to each other. We will also need to assume that $R$ and
$\rho$ are sufficiently large in comparison to $r$. This is summarized in
Assumption \ref{as:cones-setup}. The idea behind it is demonstrated in Figure
\ref{fig:M-assumptions}, which might provide some intuition.

\begin{assumption}
\label{as:cones-setup} Let $\mathbf{m}>1$ and let $\boldsymbol{\gamma}%
_{0}=(\mathbf{a}_{0},\mathbf{b}_{0},\mathbf{c}_{0})\in\mathbb{R}^{3},$
$\boldsymbol{\gamma}_{1}=(\mathbf{a}_{1},\mathbf{b}_{1},\mathbf{c}_{1}%
)\in\mathbb{R}^{3}$ satisfy $\mathbf{a}_{m}>0,$ $\mathbf{b}_{m},\mathbf{c}%
_{m}<0$ for $m=1,2.$ Let us also define a set $M\subset\mathbf{B}$ as%
\begin{equation}
M:=\overline{B}_{u}(0,r)\times\overline{B}_{s}(0,r)\times\overline{B}%
_{c}.\label{eq:Mi-set}%
\end{equation}
We assume that for any horizontal disc $\mathbf{h}$ in a ch-set with cones
$(M,\boldsymbol{\gamma}_{1})$ and for any $i\in I$ there exists $\left(
\iota,\kappa\right)  \in J$ such that $\mathbf{h}(B_{u}(0,r))\subset
\mathrm{dom}(\tilde{\eta}_{\iota}\circ\tilde{\eta}_{i}^{-1}).$ In addition we
assume that for any $q^{\ast},q^{\ast\ast}$ in $\mathrm{dom}(\tilde{\eta
}_{\iota}\circ\tilde{\eta}_{i}^{-1})$ such that $Q_{\boldsymbol{\gamma}_{1}%
}(q^{\ast}-q^{\ast\ast})>0$ we have
\begin{equation}
Q_{\boldsymbol{\gamma}_{0}}(\tilde{\eta}_{\iota}\circ\tilde{\eta}_{i}%
^{-1}(q^{\ast})-\tilde{\eta}_{\iota}\circ\tilde{\eta}_{i}^{-1}(q^{\ast\ast
}))>\mathbf{m}Q_{\boldsymbol{\gamma}_{1}}(q^{\ast}-q^{\ast\ast}%
),\label{eq:from-g1-to-g0}%
\end{equation}
and
\begin{equation}
\mathbf{h}^{\prime}:=\tilde{\eta}_{\iota}\circ\tilde{\eta}_{i}^{-1}%
\circ\mathbf{h}_{|B_{u}(0,r)}\quad\text{is a horizontal disc in }%
(M_{\iota,\kappa},\boldsymbol{\gamma}_{0}).\label{eq:h'-hor-disc}%
\end{equation}

\end{assumption}

Assumption \ref{as:cones-setup} ensures that for $\mathbf{h}$ in some local
coordinates $\tilde{\eta}_{i}$ we can change to coordinates $\tilde{\eta
}_{\iota}$ so that $\mathbf{h}^{\prime}:=\tilde{\eta}_{\iota}\circ\tilde{\eta
}_{i}^{-1}\circ\mathbf{h}$ lies close to the middle of the set $M$. Assumption
\ref{as:cones-setup} is also discussed in Section \ref{sec:local-maps}, where
conditions which imply it are given.

\begin{remark}
Above we use bold font for $\boldsymbol{\gamma}_{i}=(\mathbf{a}_{i}%
,\mathbf{b}_{i},\mathbf{c}_{i})$, $i=0,1$ to emphasize that these are fixed
constants, and to distinguish them from other $\gamma=(a,b,c)$ in our proofs.
\end{remark}

\begin{definition}
\label{def:cone-conditions}If for any $(i,k)\in J$ there exists a sequence of
ch-sets with cones $(N_{1},\gamma_{1}),\ldots,(N_{n},\gamma_{n})$ ($n$ can
depend on $(i,k)$) and a sequence $i_{0}=i,i_{1},\ldots,i_{n}\in I$ such that%
\begin{equation}
M_{i,k}=:N_{0}\overset{f_{i_{1}i_{0}}}{\Longrightarrow}N_{1}\overset
{f_{i_{2}i_{1}}}{\Longrightarrow}N_{2}\overset{f_{i_{3}i_{2}}}{\Longrightarrow
}\ldots\overset{f_{i_{n}i_{n-1}}}{\Longrightarrow}N_{n}\overset{\mathrm{id}%
}{\Longrightarrow}M,\label{eq:cover-sequence}%
\end{equation}
then we say that $f$ \emph{satisfies covering conditions}.

If in addition for any $q_{1},q_{2}\in N_{l-1},$ $q_{1}\neq q_{2},$%
\begin{equation}
Q_{\gamma_{l+1}}(f_{i_{l+1}i_{l}}(q_{1})-f_{i_{l+1}i_{l}}(q_{2}))>Q_{\gamma
_{l}}(q_{1}-q_{2})\label{eq:cc-local-iterates}%
\end{equation}
for $l=0,\ldots,n-1,$ and for $\gamma_{n}=(a,b,c)$ we have%
\begin{equation}
\mathbf{a}_{1}>a,\quad\frac{\mathbf{b}_{1}}{\mathbf{a}_{1}}>\frac{b}{a}%
,\quad\frac{\mathbf{c}_{1}}{\mathbf{a}_{1}}>\frac{c}{a},\label{eq:to-gamma0}%
\end{equation}
then we say that $f$ \emph{satisfies cone conditions}.
\end{definition}

\begin{figure}[ptb]
\begin{center}
\includegraphics[
width=5in
]{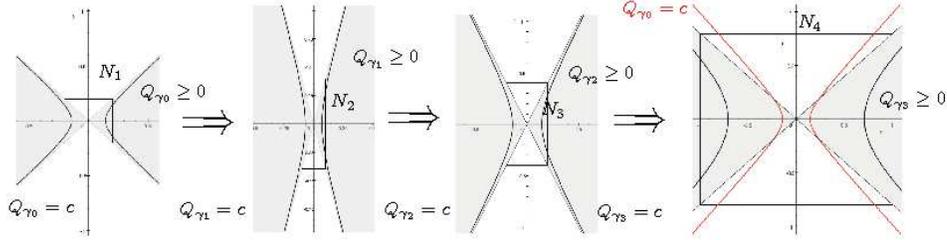}
\caption{(see Example \ref{ex:cone-cond-def}) For the first iterates of the
map the ch-sets and cones are contracted in the $x$ direction. After a number
of steps the expansion in $x$ starts to dominate. Note that the coordinate
$\theta$ is expanding. Since expansion in $x$ is stronger than expansion in
$\theta$ though, the cones eventually become more flat and their level sets
$Q_{\gamma_{i}}=c$ are pulled away from the origin.}%
\label{fig:example-cone}%
\end{center}
\end{figure}

\begin{example}
\label{ex:cone-cond-def}This example stands behind the pictures from Figure
\ref{fig:example-cone}. Consider $u=c=1$ and $s=0.$ Assume that $f_{i_{1}%
i_{0}}=(A_{ij}^{1})_{i,j=1,2}=diag(\frac{1}{2},2)$ $f_{i_{2}i_{1}}=(A_{ij}%
^{2})_{i,j=1,2}=diag(2,1),$ $f_{i_{2}i_{3}}=(A_{ij}^{3})_{i,j=1,2}=diag(5,2).$
Let $\boldsymbol{\gamma}_{0}=(1,-1)$ and $\boldsymbol{\gamma}_{1}=(\frac{1}%
{4},-\frac{3}{8}).$ We take ch-sets with cones $(N_{l}((0,0),r_{u}^{l}%
,r_{c}^{l}),\gamma_{l}),$ for $l=0,1,2,3$ with
\begin{eqnarray*}
r_{u}^{0}  &  = &r_{c}^{0}=r,\\
r_{u}^{l}  &  = &r_{u}^{l-1}A_{11}^{l}-\varepsilon,\\
r_{c}^{l}  &  = &r_{c}^{l-1}A_{22}^{l}+\varepsilon,
\end{eqnarray*}
$\gamma_{0}=\boldsymbol{\gamma}_{0},$ $\gamma_{1}=(4\delta,-\frac{1}{4}%
\delta^{-1}),$ $\gamma_{2}=(1\delta^{2},-\frac{1}{4}\delta^{-2}),$ $\gamma
_{3}=(\frac{1}{25}\delta^{3},-\frac{1}{16}\delta^{-3}),$ with $\delta
=1+\varepsilon.$ For sufficiently small $r$ and $\varepsilon$ we will have
(\ref{eq:cover-sequence}) and (\ref{eq:cc-local-iterates}). For sufficiently
small $\varepsilon$ we also have (\ref{eq:to-gamma0}). Assume now that
$\tilde{\eta}_{\iota}\circ\tilde{\eta}_{i_{3}}^{-1}=diag(1,1+\frac{1}{4}).$
This $\tilde{\eta}_{\iota}\circ\tilde{\eta}_{i_{3}}^{-1}$ is taken just as a
hypothetical example, in order to show that even when a switch to new
coordinates involves an expansion in the central coordinate the Assumption
\ref{as:cones-setup} can easily be satisfied. We have
\begin{eqnarray*}
Q_{\boldsymbol{\gamma}_{0}}(\tilde{\eta}_{\iota}\circ\tilde{\eta}_{i_{3}}%
^{-1}(x,y))  &  = & x^{2}-\frac{5}{4}\theta^{2}\\
&  = & 4\left(  \frac{1}{4}x^{2}-\frac{5}{16}\theta^{2}\right) \\
&  \geq & 4\left(  \frac{1}{4}x^{2}-\frac{3}{8}\theta^{2}\right) \\
&  = & 4Q_{\boldsymbol{\gamma}_{1}}((x,y))
\end{eqnarray*}
which means that (\ref{eq:from-g1-to-g0}) holds for $\mathbf{m}<4$.
\end{example}

We now introduce a notation $U\subset D_{\phi}$ for a set
\begin{equation}
U:=B_{u}(0,r)\times B_{s}(0,r)\times\Lambda. \label{eq:U-set-bound}%
\end{equation}
The set $U$ will be the region in which we will construct an invariant
manifold of points, which stay within the set $D_{\phi}$ for forward
iterations of the map $f_{\phi}.$

\begin{figure}[ptb]
\begin{center}
\includegraphics[
width=4.8in
]{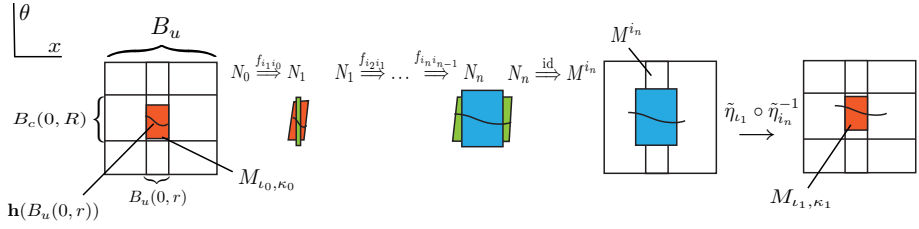}
\end{center}
\caption{The sequence of covering relations from Definition
\ref{def:cone-conditions}, together with the sets $M_{\iota_{0},\kappa_{0}}$
and $M_{\iota_{1},\kappa_{1}}$, which are the first step of the inductive
construction from the proof of Theorem \ref{th:nhims-from-cc}.}%
\label{fig:covering-iteration}%
\end{figure}

\subsection{Existence of a normally hyperbolic manifold - Main
result\label{sec:main-res}}

In this section we use the assumptions from Section
\ref{sec:cover-cc-norm-hyp} to obtain the existence of a normally hyperbolic
invariant manifold inside of the set $U$ defined in (\ref{eq:U-set-bound}). We
start with a construction of the center-stable manifold. This is given in
Theorem \ref{th:nhims-from-cc}. The existence of an center-unstable manifold
follows from mirror arguments for the inverse map. The normally hyperbolic
manifold is obtained by intersecting the center-stable and center-unstable
manifolds. This is done in Theorem \ref{th:main}.

\begin{theorem}
\label{th:nhims-from-cc}If $f$ satisfies cone conditions then there exists a
continuous monomorphism $V:B_{s}(0,r)\times\Lambda\rightarrow U$ such that

\begin{enumerate}
\item \label{it:th-nhims-cc1}$\pi_{y}V(y,\lambda)=y$, $\pi_{\lambda
}V(y,\lambda)=\lambda,$

\item \label{it:th-nhims-cc2}for any $(y,\lambda)\in B_{s}(0,r)\times\Lambda$
and any $n\in\mathbb{N}$
\[
f_{\phi}^{n}(V(y,\lambda))\in D_{\phi}.
\]

\item \label{it:th-nhims-cc3}for any $q\in U$ such that $f_{\phi}^{n}(q)\in
D_{\phi}$ for all $n\in\mathbb{N,}$ there exists a $(y,\lambda)\in
B_{s}(0,r)\times\Lambda$ such that $q=V(y,\lambda),$

\item \label{it:th-nhims-cc4}if $\lambda^{\ast},\lambda^{\ast\ast}\in\eta
_{i}^{-1}(B_{c}(\eta_{i}(\boldsymbol{\lambda}_{k}),\rho))$ for some $(i,k)\in
J$ then for any $y^{\ast},y^{\ast\ast}\in\overline{B}_{s}(0,r)$ such that
$(\lambda^{\ast},y^{\ast})\neq(\lambda^{\ast\ast},y^{\ast\ast})$%
\begin{equation}
Q_{\gamma_{0}}\left(  \tilde{\eta}_{i}\circ V(y^{\ast},\lambda^{\ast}%
)-\tilde{\eta}_{i}\circ V(y^{\ast\ast},\lambda^{\ast\ast})\right)
<0.\label{eq:V-vert-cones}%
\end{equation}

\end{enumerate}
\end{theorem}

\begin{proof}
We take any $y_{0}\in B_{s}(0,r),$ $\lambda_{0}\in\Lambda$ and $(\iota
_{0},\kappa_{0})\in J$ such that $\lambda_{0}\in\eta_{\iota_{0}}^{-1}%
(B_{c}(\boldsymbol{\lambda}_{\kappa_{0}},\rho))$ and define a horizontal disc
$\mathbf{h}_{0}$ in $M_{\iota_{0},\kappa_{0}}$ as%
\[
\mathbf{h}_{0}(x):=(x,y_{0},\eta_{\iota_{0}}(\lambda_{0})).
\]
Since $f$ satisfies cone conditions, using assumption (\ref{eq:cover-sequence}%
) and applying inductively Lemma \ref{lem:hor-disc-image} gives us the
existence of indexes $i_{1},\ldots,i_{n_{1}}\in I$ and of a horizontal disc
$\mathbf{h}_{1}$ in $(M,\gamma_{n_{1}})$ such that
\begin{eqnarray*}
\mathbf{h}_{1}(\overline{B}_{u}) & = &  \{\tilde{\eta}_{i_{n_{1}}}\circ f_{\phi
}^{n_{1}}\circ\tilde{\eta}_{\iota_{0}}^{-1}(\mathbf{h}_{0}(x))\in M:x\in
B_{u}(0,r),\text{ and}\\
& &  \quad\tilde{\eta}_{i_{l}}\circ f_{\phi}^{l}\circ\tilde{\eta}_{\iota_{0}%
}^{-1}(\mathbf{h}_{1}(x))\in N_{l}\text{ for }l=1,\ldots,n_{1}\}.
\end{eqnarray*}
By (\ref{eq:to-gamma0})
\[
Q_{\mathbf{\gamma}_{1}}(\mathbf{h}_{1}(x^{\ast})-\mathbf{h}_{1}(x^{\ast\ast
}))>Q_{\gamma_{n_{1}}}(\mathbf{h}_{1}(x^{\ast})-\mathbf{h}_{1}(x^{\ast\ast
}))>0,
\]
which means that $\mathbf{h}_{1}$ is a horizontal disc in
$(M,\boldsymbol{\gamma}_{1}).$ From (\ref{eq:from-g1-to-g0}) and
(\ref{eq:h'-hor-disc}) we know that there exists $(\iota_{1},\kappa_{1})\in J$
such that $\mathbf{h}^{\prime}:=\tilde{\eta}_{\iota_{1}}\circ\tilde{\eta
}_{i_{n_{1}}}^{-1}\circ\mathbf{h}_{1}$ is a horizontal disc in $(M_{\iota
_{1},\kappa_{1}},\boldsymbol{\gamma}_{0}).$ This in particular means that for
$\mathbf{f}_{1}:=\tilde{\eta}_{\iota_{1}}\circ f_{\phi}^{n_{1}}\circ
\tilde{\eta}_{\iota_{0}}^{-1},$ there exists an $x\in B_{u}(0,r)$ for which
$\mathbf{f}_{1}(\mathbf{h}_{1}(x))\cap M_{\iota_{1},\kappa_{1}}\neq\emptyset.$
By (\ref{eq:cc-local-iterates}) and (\ref{eq:from-g1-to-g0}), for any
$x^{\ast}\neq x^{\ast\ast}$ such that $\mathbf{h}_{1}(x^{\ast}),\mathbf{h}%
_{1}(x^{\ast\ast})\in\mathrm{dom}(\mathbf{f}_{1})$
\begin{equation}
Q_{\gamma_{0}}\left(  \mathbf{f}_{1}(\mathbf{h}_{1}(x^{\ast}))-\mathbf{f}%
_{1}(\mathbf{h}_{1}(x^{\ast\ast}))\right)  >\mathbf{m}Q_{\gamma_{0}}\left(
\mathbf{h}_{1}(x^{\ast})-\mathbf{h}_{1}(x^{\ast\ast})\right)
>0.\label{eq:cc-inductive}%
\end{equation}

Repeating the above procedure inductively (starting the second step with the
horizontal disc $\mathbf{h}^{\prime}$ and local coordinates given by
$\tilde{\eta}_{\iota_{1}}$) we obtain a sequence of points $x_{s}\in
B_{u}(0,r)$ and indexes $(\iota_{s},\kappa_{s})$ for $s\in\mathbb{N}$ such
that for
\[
\mathbf{f}_{s}:=\tilde{\eta}_{\iota_{s}}\circ f_{\phi}^{n_{s}+\ldots+n_{1}%
}\circ\tilde{\eta}_{\iota_{0}}^{-1}%
\]
we have
\[
\mathbf{f}_{w}(\mathbf{h}_{1}(x_{s}))\in M_{\iota_{w},\kappa_{w}}\quad
\quad\text{for }w\leq s.
\]
Since $\overline{B}_{u}(0,r)$ is compact, there exists an $x_{0}=x_{0}%
(y_{0},\lambda_{0})\in B_{u}(0,r)$ such that $\tilde{\eta}_{\iota_{s}}%
^{-1}\circ\mathbf{f}_{s}(\mathbf{h}_{1}(x_{0}))\in U$ for all $s\in
\mathbb{N.}$ We define $V(y_{0},\lambda_{0}):=\tilde{\eta}_{\iota_{0}}%
^{-1}(x_{0}(y_{0},\lambda_{0}),y_{0},\lambda_{0}).$ To see that $V$ is
properly defined suppose that we have two points $x_{0}^{\ast}\neq x_{0}%
^{\ast\ast}$ such that $\tilde{\eta}_{\iota_{s}}\circ\mathbf{f}(\mathbf{h}%
_{1}(x_{0}^{\ast})),\tilde{\eta}_{\iota_{s}}\circ\mathbf{f}(\mathbf{h}%
_{1}(x_{0}^{\ast\ast}))\in U$ for all $s\in\mathbb{N.}$ Then by
(\ref{eq:cc-inductive}) we obtain%
\begin{eqnarray}
Q_{\gamma_{0}}\left(  \mathbf{f}_{s}(\mathbf{h}_{1}(x_{0}^{\ast}%
))-\mathbf{f}_{s}(\mathbf{h}_{1}(x_{0}^{\ast\ast}))\right)   &  > &\mathbf{m}%
Q_{\gamma_{0}}\left(  \mathbf{f}_{s-1}(\mathbf{h}_{1}(x_{0}^{\ast
}))-\mathbf{f}_{s-1}(\mathbf{h}_{1}(x_{0}^{\ast\ast}))\right)  \nonumber\\
&  >&\ldots\label{eq:h-pulled-ap}\\
&  >&\mathbf{m}^{s}Q_{\gamma_{0}}\left(  \mathbf{h}_{1}(x_{0}^{\ast
})-\mathbf{h}_{1}(x_{0}^{\ast\ast})\right)  \nonumber\\
&  >&0.\nonumber
\end{eqnarray}
Since $\mathbf{m}>1,$ (\ref{eq:h-pulled-ap}) implies in particular that
\[
\left\Vert \pi_{x}\left(  \mathbf{f}_{s}(\mathbf{h}_{1}(x_{0}^{1}%
))-\mathbf{f}_{s}(\mathbf{h}_{1}(x_{0}^{2}))\right)  \right\Vert
\rightarrow\infty\quad\text{as \quad}s\rightarrow\infty.
\]
This is impossible since $\mathbf{f}_{s}(\mathbf{h}_{1}(x_{0}^{w}))$ is in
$M_{\iota_{s},\kappa_{s}},$ which is a subset of $\mathbf{B,}$ which is a bounded.

We now need to show (\ref{eq:V-vert-cones}). Suppose that $V(y^{\ast}%
,\lambda^{\ast}),V(y^{\ast\ast},\lambda^{\ast\ast})\in M_{i,k}$ and
$Q_{\gamma_{0}}\left(  \tilde{\eta}_{i}\circ V(y^{\ast},\lambda^{\ast}%
)-\tilde{\eta}_{i}\circ V(y^{\ast\ast},\lambda^{\ast\ast})\right)  \geq0.$
Applying estimates analogous to (\ref{eq:h-pulled-ap}) we obtain a contradiction.

Continuity of $V$ will follow from the fact that
\begin{equation}
Q_{\gamma_{0}}\left(  \tilde{\eta}_{\iota_{0}}\circ V(y^{\ast},\lambda^{\ast
})-\tilde{\eta}_{\iota_{0}}\circ V(y^{\ast\ast},\lambda^{\ast\ast})\right)
<0.\label{ea:cont-V-cone}%
\end{equation}
Since $\boldsymbol{\gamma}_{0}=(\mathbf{a}_{0},\mathbf{b}_{0},\mathbf{c}_{0})$
with $\mathbf{a}_{0}>0$ and $\mathbf{b}_{0},\mathbf{c}_{0}<0$
(\ref{ea:cont-V-cone}) gives
\begin{eqnarray*}
0 &  >&Q_{\gamma_{0}}\left(  \tilde{\eta}_{\iota_{0}}\circ V(y^{\ast}%
,\lambda^{\ast})-\tilde{\eta}_{\iota_{0}}\circ V(y^{\ast\ast},\lambda
^{\ast\ast})\right)  \\
&  = &\mathbf{a}_{0}\left\Vert \pi_{x}V(y^{\ast},\lambda^{\ast})-\pi
_{x}V(y^{\ast\ast},\lambda^{\ast\ast})\right\Vert ^{2}+\mathbf{b}%
_{0}\left\Vert y^{\ast}-y^{\ast\ast}\right\Vert ^{2} \\ 
& & +\mathbf{c}_{0}\left\Vert
\eta_{\iota_{0}}(\lambda^{\ast})-\eta_{\iota_{0}}(\lambda^{\ast\ast
})\right\Vert ^{2},
\end{eqnarray*}
and therefore%
\begin{eqnarray*}
\mathbf{a}_{0}\left\Vert \pi_{x}V(y^{\ast},\lambda^{\ast})-\pi_{x}%
V(y^{\ast\ast},\lambda^{\ast\ast})\right\Vert ^{2 } \\
<\min(-\mathbf{b}%
_{0},-\mathbf{c}_{0})\left\Vert \left(  y^{\ast},\eta_{\iota_{0}}%
(\lambda^{\ast})\right)  -\left(  y^{\ast\ast},\eta_{\iota_{0}}(\lambda
^{\ast\ast})\right)  \right\Vert ^{2}.
\end{eqnarray*}

\end{proof}

Now we move to proving the existence of the normally hyperbolic invariant
manifold. First we need a definition.

\begin{definition}
We say that $f$ \emph{satisfies backward cone conditions} if $f^{-1}$
satisfies cone conditions, with reversed roles of $x$ and $y$ coordinates.
\end{definition}

We assume that for $f$ Assumption \ref{as:cones-setup} holds with
$\boldsymbol{\gamma}_{0}=\boldsymbol{\gamma}_{0}^{\mathrm{forw}}.$ We assume
also that for $f^{-1}$ Assumption \ref{as:cones-setup} holds with
$\boldsymbol{\gamma}_{0}=\boldsymbol{\gamma}_{0}^{\mathrm{back}}$ (with
reversed roles of the $x$ and $y$ coordinates).

\begin{theorem}
\label{th:main}(Main Theorem) Assume that $f$ satisfies cone conditions for
$\boldsymbol{\gamma}_{0}^{\mathrm{forw}}=\left(  \mathbf{a}_{0}^{\mathrm{f}%
},\mathbf{b}_{0}^{\mathrm{f}},\mathbf{c}_{0}^{\mathrm{f}}\right)  $ and
backward cone conditions with $\boldsymbol{\gamma}_{0}^{\mathrm{back}}=\left(
\mathbf{a}_{0}^{\mathrm{b}},\mathbf{b}_{0}^{\mathrm{b}},\mathbf{c}%
_{0}^{\mathrm{b}}\right)  .$ If%
\begin{equation}
\left\vert \mathbf{a}_{0}^{\mathrm{f}}\right\vert >\left\vert \mathbf{a}%
_{0}^{\mathrm{b}}\right\vert \text{\quad and\quad}\left\vert \mathbf{b}%
_{0}^{\mathrm{f}}\right\vert <\left\vert \mathbf{b}_{0}^{\mathrm{b}%
}\right\vert \label{eq:ab-ineq}%
\end{equation}
then there exist continuous monomorphisms $W^{s}:B_{s}(0,r)\times
\Lambda\rightarrow U,$ $W^{u}:B_{u}(0,r)\times\Lambda\rightarrow U$ and
$\chi:\Lambda\rightarrow U,$ such that
\begin{equation}
\pi_{y,\lambda}W^{s}(y,\lambda)=(y,\lambda),\quad\pi_{x,\lambda}%
W^{u}(y,\lambda)=(x,\lambda),\quad\pi_{\lambda}\chi(\lambda)=\lambda
,\label{eq:W-projections}%
\end{equation}
and $\Lambda_{\phi}:=\chi(\Lambda)$ is an invariant manifold for $f_{\phi},$
with stable manifold $W^{s}(B_{s}(0,r)\times\Lambda)$ and unstable manifold
$W^{u}(B_{u}(0,r)\times\Lambda).$
\end{theorem}

\begin{proof}
Since $f$ satisfies cone conditions, applying Theorem \ref{th:nhims-from-cc}
we obtain $W^{s}(y,\lambda)$ as $V$. Since $f$ satisfies backward cone
conditions, once again from Theorem \ref{th:nhims-from-cc} for $f^{-1}$ we
also obtain $W^{u}(x,\lambda)$ as function $V$. From point
\ref{it:th-nhims-cc1} in Theorem \ref{th:nhims-from-cc} it follows that
(\ref{eq:W-projections}) holds for $W^{s}$ and $W^{u}$.

We shall show that for any $\lambda\in\Lambda$ the sets $W^{s}(B_{s}%
(0,r),\lambda)$ and $W^{u}(B_{u}(0,r),\lambda)$ intersect. Let us define
$F:B_{u}(0,r)\times B_{s}(0,r)\rightarrow B_{u}(0,r)\times B_{s}(0,r)$ as%
\[
F(x,y):=\left(  \pi_{x}W^{s}(y,\lambda),\pi_{y}W^{u}(x,\lambda)\right)  .
\]
Since $F$ is continuous, from the Brouwer fixed point theorem follows that
there exists an $(x_{0},y_{0})$ such that $F(x_{0},y_{0})=\left(  x_{0}%
,y_{0}\right)  .$ By (\ref{eq:W-projections}) this means that%
\[
W^{s}(y_{0},\lambda)=\left(  \pi_{x}W^{s}(y_{0},\lambda),y_{0},\lambda\right)
=\left(  x_{0},\pi_{y}W^{u}(x_{0},\lambda),\lambda\right)  =W^{u}%
(x_{0},\lambda).
\]

Now we shall show that for any given $\lambda\in\Lambda$ there exists only a
single point of such intersection. Suppose that for some $\lambda\in\Lambda$
there exist $\left(  x^{\ast},y^{\ast}\right)  ,\left(  x^{\ast\ast}%
,y^{\ast\ast}\right)  \in B_{u}(0,r)\times B_{s}(0,r)$, $\left(  x^{\ast
},y^{\ast}\right)  \neq\left(  x^{\ast\ast},y^{\ast\ast}\right)  $ such that
\[
W^{s}(y^{\ast},\lambda)=W^{u}(x^{\ast},\lambda)\quad\text{and\quad}%
W^{s}(y^{\ast\ast},\lambda)=W^{u}(x^{\ast\ast},\lambda).
\]
From (\ref{eq:W-projections}) we have $W^{s}(y_{m},\lambda)=W^{u}%
(x_{m},\lambda)=(x_{m},y_{m},\lambda)$ for $m=1,2.$ From point 4. in Theorem
\ref{th:nhims-from-cc} follows that%
\begin{eqnarray*}
Q_{\gamma_{0}^{\mathrm{forw}}}\left(  \tilde{\eta}_{i}\circ W^{s}(y^{\ast
},\lambda)-\tilde{\eta}_{i}\circ W^{s}(y^{\ast\ast},\lambda)\right)  =\\
\qquad Q_{\gamma_{0}^{\mathrm{forw}}}\left(  (x^{\ast},y^{\ast},\eta_{i}%
(\lambda))-(x^{\ast\ast},y^{\ast\ast},\eta_{i}(\lambda))\right)  <0,
\end{eqnarray*}%
\begin{eqnarray*}
Q_{\gamma_{0}^{\mathrm{back}}}\left(  \tilde{\eta}_{i}\circ W^{u}(x^{\ast
},\lambda)-\tilde{\eta}_{i}\circ W^{u}(x^{\ast\ast},\lambda)\right)  =\\
\qquad Q_{\gamma_{0}^{\mathrm{back}}}\left(  (x^{\ast},y^{\ast},\eta_{i}%
(\lambda))-(x^{\ast\ast},y^{\ast\ast},\eta_{i}(\lambda))\right)  <0.
\end{eqnarray*}
which implies that
\begin{eqnarray}
a_{0}^{\mathrm{f}}\left\Vert x^{\ast}-x^{\ast\ast}\right\Vert ^{2}%
+b_{0}^{\mathrm{f}}\left\Vert y^{\ast}-y^{\ast\ast}\right\Vert ^{2} &
< & 0,\label{eq:temp-contr-1}\\
a_{0}^{\mathrm{b}}\left\Vert x^{\ast}-x^{\ast\ast}\right\Vert ^{2}%
+b_{0}^{\mathrm{b}}\left\Vert y^{\ast}-y^{\ast\ast}\right\Vert ^{2} &
< & 0.\label{eq:temp-contr-2}%
\end{eqnarray}
From (\ref{eq:ab-ineq}) and (\ref{eq:temp-contr-2}) (keeping in mind that
$a_{0}^{\mathrm{f}}>0,$ $b_{0}^{\mathrm{f}}<0$ and that $a_{0}^{\mathrm{b}%
}<0,$ $b_{0}^{\mathrm{b}}>0$ due to the reversion of the roles of $x$ and $y$
for the inverse map) follows that
\[
a_{0}^{\mathrm{f}}\left\Vert x^{\ast}-x^{\ast\ast}\right\Vert ^{2}%
>-a_{0}^{\mathrm{b}}\left\Vert x^{\ast}-x^{\ast\ast}\right\Vert ^{2}%
>b_{0}^{\mathrm{b}}\left\Vert y^{\ast}-y^{\ast\ast}\right\Vert ^{2}%
>-b_{0}^{\mathrm{f}}\left\Vert y^{\ast}-y^{\ast\ast}\right\Vert ^{2},
\]
which contradicts (\ref{eq:temp-contr-1}).

We now define $\chi(\lambda):=(x_{0},y_{0},\lambda)$ for $x_{0},$ $y_{0}$ such
that $W^{s}(y_{0},\lambda)=W^{u}(x_{0},\lambda).$ By above arguments we know
that $\chi$ is a properly defined function. We need to show that this function
is continuous. Let us take any $\lambda^{\ast},\lambda^{\ast\ast}\in\eta
_{i}^{-1}(B_{c}(\eta_{i}(\boldsymbol{\lambda}_{k}),\rho))$ for some $(i,k)\in
J.$ From point 4 in Theorem \ref{th:nhims-from-cc} follows that%
\begin{eqnarray}
Q_{\gamma_{0}^{\mathrm{forw}}}\left(  \tilde{\eta}_{i}\circ\chi(\lambda^{\ast
})-\tilde{\eta}_{i}\circ\chi(\lambda^{\ast\ast})\right)   &
<&0,\label{eq:cc-for-chi-2}\\
Q_{\gamma_{0}^{\mathrm{back}}}\left(  \tilde{\eta}_{i}\circ\chi(\lambda^{\ast
})-\tilde{\eta}_{i}\circ\chi(\lambda^{\ast\ast})\right)   &  <&0.\nonumber
\end{eqnarray}
Let us adopt notations $\tilde{\eta}_{i}\circ\chi(\lambda^{\ast})=\left(
x^{\ast},y^{\ast},\theta^{\ast}\right)  $ and $\tilde{\eta}_{i}\circ
\chi(\lambda^{\ast\ast})=\left(  x^{\ast\ast},y^{\ast\ast},\theta^{\ast\ast
}\right)  .$ Note that from the construction of $\chi$ follows that $\eta
_{i}(\lambda^{\ast})=\theta^{\ast}$ and $\eta_{i}(\lambda^{\ast\ast}%
)=\theta^{\ast\ast}.$ From (\ref{eq:cc-for-chi-2}) it follows that
\begin{eqnarray}
&&  \left(  a_{0}^{\mathrm{f}}+a_{0}^{\mathrm{b}}\right)  \left\Vert x^{\ast
}-x^{\ast\ast}\right\Vert ^{2}+\left(  b_{0}^{\mathrm{f}}+b_{0}^{\mathrm{b}%
}\right)  \left\Vert y^{\ast}-y^{\ast\ast}\right\Vert ^{2}%
\label{eq:cc-for-chi}\\
&&  <-\left(  c_{0}^{\mathrm{f}}+c_{0}^{\mathrm{b}}\right)  \left\Vert
\theta^{\ast}-\theta^{\ast\ast}\right\Vert ^{2}\nonumber\\
&&  =-\left(  c_{0}^{\mathrm{f}}+c_{0}^{\mathrm{b}}\right)  \left\Vert \eta
_{i}(\lambda^{\ast})-\eta_{i}(\lambda^{\ast\ast})\right\Vert ^{2}\nonumber
\end{eqnarray}
From (\ref{eq:ab-ineq}) it follows that $a_{0}^{\mathrm{f}}+a_{0}^{\mathrm{b}%
}=\left\vert a_{0}^{\mathrm{f}}\right\vert -\left\vert a_{0}^{\mathrm{b}%
}\right\vert >0$ and $b_{0}^{\mathrm{f}}+b_{0}^{\mathrm{b}}=-\left\vert
b_{0}^{\mathrm{f}}\right\vert +\left\vert b_{0}^{\mathrm{b}}\right\vert >0.$
By the fact that $\eta_{i}$ is continuous and the fact that $c_{0}%
^{\mathrm{f}}<0$ and $c_{0}^{\mathrm{b}}<0$, from (\ref{eq:cc-for-chi})
follows the continuity of $\chi.$

We will now show that for any $p\in W^{s}(B_{s}(0,r)\times\Lambda),$ $f_{\phi
}^{n}(p)$ converges to $\chi(\Lambda)$ as $n$ goes to infinity. Let us
consider the limit set of the point $p$%
\[
\omega(f_{\phi},p)=\{q|\lim_{k\rightarrow\infty}f_{\phi}^{n_{k}}(p)=q\text{
for some }n_{k}\rightarrow\infty\}.
\]
If we can show that $\omega(f_{\phi},p)$ is contained in $W^{u}\cap W^{s}%
=\chi(\Lambda),$ then this will conclude our proof. We take any $q=\lim
_{k\rightarrow\infty}f_{\phi}^{n_{k}}(p)$ from $\omega(f_{\phi},p).$ By
continuity of $W^{s}$ we know that $q\in W^{s}.$ Suppose now that $q\notin
W^{u}.$ This would mean that there exists an $n>0$ for which $f_{\phi}%
^{-n}(q)\notin B_{u}(0,r)\times B_{s}(0,r)\times\Lambda.$ Since%
\[
\lim_{k\rightarrow\infty}f_{\phi}^{n_{k}-n}(p)=f_{\phi}^{-n}(q),
\]
we have that $f_{\phi}^{-n}(q)\in\omega(f_{\phi},p),$ but this contradicts the
fact that $\omega(f_{\phi},p)\subset B_{u}(0,r)\times B_{s}(0,r)\times
\Lambda.$

Showing that all backward iterations of points in $W^{u}(B_{u}(0,r)\times
\Lambda)$ converge to $\chi(\Lambda)$ is analogous.
\end{proof}

\begin{remark}
Let us note that during the course of the proof of Theorem \ref{th:main} we
have established more than just continuity of $W^{u},$ $W^{s}$ and $\chi.$
From our construction we know that for $i\in I$%
\begin{eqnarray*}
\tilde{\eta}_{i}\circ W^{u}(x,\eta_{i}^{-1}(\theta))  &  = &\left(  x,w_{i}%
^{u}(x,\theta),\theta\right)  ,\\
\tilde{\eta}_{i}\circ W^{s}(y,\eta_{i}^{-1}(\theta))  &  = &\left(  w_{i}%
^{s}(y,\theta),y,\theta\right)  ,\\
\tilde{\eta}_{i}\circ\chi(\eta_{i}^{-1}(\theta))  &  = & \left(  \varkappa
_{i}(\theta),\theta\right)  ,
\end{eqnarray*}
for continuous $w_{i}^{u}:B_{u}(0,r)\times B_{c}\rightarrow B_{s}(0,r),$
$w_{i}^{s}:B_{s}(0,r)\times B_{c}\rightarrow B_{u}(0,r)$ and $\varkappa
_{i}:B_{c}\rightarrow B_{u}(0,r)\times B_{s}(0,r).$ The inequality
(\ref{eq:V-vert-cones}) from Theorem \ref{th:nhims-from-cc} can be used to
obtain explicit Lipschitz bounds for functions $w_{i}^{u},$ $w_{i}^{s}.$ Also
estimates (\ref{eq:cc-for-chi}) can be used to obtain Lipschitz bounds for
$\varkappa_{i}.$ This means that we can get Lipschitz estimates for the
invariant manifold $\chi(\Lambda)$ together with Lipschitz estimates for its
stable and unstable manifold.
\end{remark}

%% file: 04-verif.tex


\section{Verification of covering and cone conditions}

\label{sec:ver-cond}

In this section we show how covering relations and cone conditions can be
verified with the use of local bounds on derivatives. The idea is to develop a
simple automatised scheme which could be applied in computer assisted proofs.
In our approach we set up our verification so that we do not need to compute
images of large sets (which in case of rigorous numerics is always
troublesome). The scheme is based on iterates of a number of single points,
combined with estimates on derivatives around their neighbourhoods.

For any set $V\subset\mathbb{R}^{n}$ we define the interval enclosure of the
derivative of $f$ on $V$ as%
\begin{eqnarray*}
\lbrack df(V)]:= \\
\left\{  A\in\mathbb{R}^{n\times n}|A_{ij}\in\left[
\inf_{x\in V}\frac{df_{i}}{dx_{j}}(x),\sup_{x\in V}\frac{df_{i}}{dx_{j}%
}(x)\right] \text{ for all }i,j=1,\ldots,n\right\}  .
\end{eqnarray*}
Let $U_{i_{1}},U_{i_{2}}\subset\Lambda$ be such that $\mathrm{dom}%
f_{i_{2}i_{1}}$ is nonempty. Assume that for any $(c+u+s)\times(c+u+s)$
matrix
\begin{equation}
A\in\lbrack df_{i_{2}i_{1}}(\mathrm{dom}f_{i_{2}i_{1}})]
\label{eq:A-for-i1-i2}%
\end{equation}
we have the following bounds%
\begin{eqnarray}
\sup\left\{  \left\Vert A_{ij}v_{j}\right\Vert :\left\Vert v_{j}\right\Vert
=1\right\}   &  \leq &\overline{A}_{ij}\label{eq:der-bounds}\\
\inf\left\{  \left\Vert A_{ij}v_{j}\right\Vert :\left\Vert v_{j}\right\Vert
=1\right\}   &  \geq & \underline{A}_{ij},\nonumber
\end{eqnarray}
with $i,j\in\{1,2,3\}$ and $v_{1},v_{2},v_{3}$ representing the variables
$x,y,\theta$ respectively (note that $\overline{A}_{ij},\underline{A}_{ij}$
depend on the choice of $i_{2}i_{1}$). In this section we shall use the bounds
(\ref{eq:der-bounds}) for verification of covering and cone conditions.

\subsection{Verifying covering conditions\label{sec:ver-cover}}

We define a $3\times3$ matrix $T_{i_{2}i_{1}}$ as%
\[
T_{i_{2}i_{1}}:=\left(  t_{ij}\right)  _{i,j=1,\ldots,3}%
\]%
\begin{equation}%
\begin{array}
[c]{lll}%
t_{11}=\underline{A}_{11},\quad & t_{12}=-\overline{A}_{12},\quad &
t_{13}=-\overline{A}_{13},\\
t_{21}=\overline{A}_{21}, & t_{22}=\overline{A}_{22}, & t_{23}=\overline
{A}_{23},\\
t_{31}=\overline{A}_{31}, & t_{32}=\overline{A}_{32}, & t_{33}=\overline
{A}_{33}.
\end{array}
\label{eq:Tmatrix-bounds}%
\end{equation}

We will use notations $R=(r_{u},r_{s},r_{c})\in\mathbb{R}^{3}$ and for
$q=(x,y,\theta)\in\mathbb{R}^{u}\times\mathbb{R}^{s}\times\mathbb{R}^{c}$ and
write%
\[
N(q,R):=N(q,r_{u},r_{s},r_{c}).
\]
We give a lemma, which can be used in order to verify that $N_{1}%
\overset{f_{i_{2}i_{1}}}{\Longrightarrow}N_{2}$.

\begin{lemma}
\label{lem:matrix-bounds-covering}Let $\varepsilon>0$ be a small number. Let
$N_{1}=N(q_{1},R_{1})\subset\mathrm{dom}f_{i_{2}i_{1}}$ be a ch-set. If for
$R_{2}=(r_{u}^{2},r_{s}^{2},r_{c}^{2}):=T_{i_{2}i_{1}}R_{1}+(-\varepsilon
,\varepsilon,\varepsilon)$ we have $r_{u}^{2},r_{s}^{2},r_{c}^{2}>0$ and for
$q_{2}:=f_{i_{2}i_{1}}(q_{1})$%
\begin{equation}
\left\Vert \pi_{x}q_{2}\right\Vert +r_{u}^{2}\leq1,\qquad\left\Vert \pi
_{y}q_{2}\right\Vert +r_{s}^{2}\leq1,\qquad\left\Vert \pi_{\theta}%
q_{2}\right\Vert +r_{c}^{2}\leq1, \label{eq:N2-in-B}%
\end{equation}
then for $N_{2}:=N(q_{2},R_{2})$ we have $N_{1}\overset{f_{i_{2}i_{1}}%
}{\Longrightarrow}N_{2}.$
\end{lemma}

\begin{proof}
Condition (\ref{eq:covering-cond-1}) holds by the choice of $q_{2}$ and
$N_{2}.$ Let $q\in N_{1}^{-},$ then for%
\[
A:=\int_{0}^{1}Df_{i_{2}i_{1}}(q_{1}+t(q-q_{1}))dt\in\lbrack df_{i_{2}i_{1}%
}(\mathrm{dom}f_{i_{2}i_{1}})],
\]
we have estimates%
\begin{eqnarray*}
\left\Vert \pi_{x}(f_{i_{2}i_{1}}(q)-q_{2})\right\Vert  &  = &\left\Vert \pi
_{x}(f_{i_{2}i_{1}}(q)-f_{i_{2}i_{1}}(q_{1}))\right\Vert \\
&  =&\left\Vert \pi_{x}\left(  \int_{0}^{1}Df_{i_{2}i_{1}}(q_{1}+t(q-q_{1}%
))dt\cdot(q-q_{1})\right)  \right\Vert \\
&  =&\left\Vert \pi_{x}A(q-q_{1})\right\Vert \\
&  =&\left\Vert A_{11}\pi_{x}(q-q_{1})+A_{12}\pi_{y}(q-q_{1})+A_{13}\pi
_{\theta}(q-q_{1})\right\Vert \\
&  \geq&\underline{A}_{11}r_{u}^{1}-\overline{A}_{12}r_{s}^{1}-\overline
{A}_{13}r_{c}^{1}\\
&  >&r_{u}^{2},
\end{eqnarray*}
hence (\ref{eq:covering-cond-2}) holds. Analogous computations for $q\in
N_{1}$ give%
\begin{eqnarray*}
\left\Vert \pi_{y}(f_{i_{2}i_{1}}(q)-q_{2})\right\Vert  &  = &\left\Vert \pi
_{y}A(q-q_{1})\right\Vert \leq\overline{A}_{21}r_{u}^{1}+\overline{A}%
_{22}r_{s}^{1}+\overline{A}_{23}r_{c}^{1}<r_{s}^{2},\\
\left\Vert \pi_{\theta}(f_{i_{2}i_{1}}(q)-q_{2})\right\Vert  &  =& \left\Vert
\pi_{\theta}A(q-q_{1})\right\Vert \leq\overline{A}_{31}r_{u}^{1}+\overline
{A}_{32}r_{s}^{1}+\overline{A}_{33}r_{c}^{1}<r_{c}^{2},
\end{eqnarray*}
which proves (\ref{eq:covering-cond-3}). Conditions (\ref{eq:N2-in-B}) ensure
that $N_{2}\subset\mathbf{B}$.
\end{proof}

\begin{example}
We return to our Example \ref{ex:cone-cond-def}. The ch-sets from the example
follow from Lemma \ref{def:cone-conditions} as $N_{l}=N(0,R_{l})$ where
$R_{0}=(r,r)$ and $R_{l+1}=T_{i_{l+1}i_{l}}R_{l}+(-\varepsilon,\varepsilon
,\varepsilon)$ with $T_{i_{l+1}i_{l}}=\mathrm{diag}(A_{11}^{l+1},A_{22}%
^{l+1}).$
\end{example}

\begin{remark}
When the $x$ coordinate is strongly expanding, for practical reasons it might
be beneficial to set\textbf{ }$r_{u}^{2}$ significantly smaller than $\pi
_{1}T_{i_{2}i_{1}}R_{1}.$ In such case the covering $N_{1}\overset
{f_{i_{2}i_{1}}}{\Longrightarrow}N_{2}$ will still take place, but $N_{2}$
will be a smaller set. This might give better bounds for next iterations of
the map $f$ and also keep the later constructed $N_{i}$ within $\mathbf{B}$.
Without reducing $r_{u}$, in the case when $x$ is expanding, it might turn out
that the sets $N_{i}$ blow up quickly.
\end{remark}

\subsection{Verifying cone conditions\label{sec:ver-cc}}

Now we shall present some lemmas, which will show how one can obtain condition
(\ref{eq:cc-local-iterates}), from bounds on derivatives (\ref{eq:der-bounds}%
). The aim is to present a simple mechanism in which successive $\gamma_{l}$
are constructed.

Let $C=(c_{ij})_{i,j=1,...,3}$ be a $3\times3$ matrix with coefficients%
\begin{equation}%
\begin{array}
[c]{lll}%
c_{11}=\underline{A}_{11}^{2}-\sum_{k\neq1}\overline{A}_{11}\overline{A}%
_{1k}\quad & c_{12}=\sum_{k=1}^{3}\overline{A}_{21}\overline{A}_{2k}\quad &
c_{13}=\sum_{k=1}^{3}\overline{A}_{31}\overline{A}_{3k}\\
c_{21}=\underline{A}_{12}^{2}-\sum_{k\neq2}\overline{A}_{12}\overline{A}%
_{1k}\quad & c_{22}=\sum_{k=1}^{3}\overline{A}_{22}\overline{A}_{2k}\quad &
c_{23}=\sum_{k=1}^{3}\overline{A}_{32}\overline{A}_{3k}\\
c_{31}=\underline{A}_{13}^{2}-\sum_{k\neq3}\overline{A}_{13}\overline{A}%
_{1k}\quad & c_{32}=\sum_{k=1}^{3}\overline{A}_{23}\overline{A}_{2k}\quad &
c_{33}=\sum_{k=1}^{3}\overline{A}_{33}\overline{A}_{3k}%
\end{array}
\label{eq:C-def}%
\end{equation}
(note that $C$ depends on the choice of $i_{2},i_{1}$).

We start with a technical lemma

\begin{lemma}
\label{lem:matrix-bounds}Let $\gamma=(a,b,c)\in\mathbb{R}^{3}$ and let
$A:\mathbb{R}^{u+s+c}\rightarrow\mathbb{R}^{u+s+c}$ be a matrix for which the
bounds (\ref{eq:der-bounds}) hold. If $a\geq0,b\leq0,c\leq0$ then for any
$p=\left(  p_{1},p_{2},p_{3}\right)  \in\mathbb{R}^{c}\times\mathbb{R}%
^{u}\times\mathbb{R}^{s}$%
\[
Q_{\gamma}(Ap)\geq Q_{C\gamma}\left(  p\right)  .
\]

\end{lemma}

\begin{proof}
Using the estimate%
\[
\pm2\left\langle A_{ki}p_{i},A_{kj}p_{j}\right\rangle \geq-\overline{A}%
_{ki}\overline{A}_{kj}\left(  \left\Vert p_{i}\right\Vert ^{2}+\left\Vert
p_{j}\right\Vert ^{2}\right)
\]
we obtain
\begin{eqnarray*}
 Q_{\gamma}(Ap) \\
  =a\sum_{i,j=1}^{3}\left\langle A_{1i}p_{i},A_{1j}%
p_{j}\right\rangle +b\sum_{i,j=1}^{3}\left\langle A_{2i}p_{i},A_{2j}%
p_{j}\right\rangle +c\sum_{i,j=1}^{3}\left\langle A_{3i}p_{i},A_{3j}%
p_{j}\right\rangle \\
  = a\sum_{i=1}^{3}||A_{1i}p_{i}||^{2}+b\sum_{i=1}^{3}||A_{2i}p_{i}%
||^{2}+c\sum_{i=1}^{3}||A_{3i}p_{i}||^{2}\\
  \quad+2\sum_{i<j}a\left\langle A_{1i}p_{i},A_{1j}p_{j}\right\rangle
+2\sum_{i<j}b\left\langle A_{2i}p_{i},A_{2j}p_{j}\right\rangle +2\sum
_{i<j}c\left\langle A_{3i}p_{i},A_{3j}p_{j}\right\rangle \\
  \geq\left\Vert p_{1}\right\Vert ^{2}(a\underline{A}_{11}^{2}+b\overline
{A}_{21}^{2}+c\overline{A}_{31}^{2})+\left\Vert p_{2}\right\Vert ^{2}(a\underline{A}_{12}^{2}%
+b\overline{A}_{22}^{2}+c\overline{A}_{32}^{2})\\
 \quad+\left\Vert p_{3}\right\Vert ^{2}(a\underline{A}_{13}^{2}%
+b\overline{A}_{23}^{2}+c\overline{A}_{33}^{2})\\
  \quad-a\sum_{i<j}\overline{A}_{1i}\overline{A}_{1j}\left(  ||p_{i}%
||^{2}+||p_{j}||^{2}\right)  +b\sum_{i<j}\overline{A}_{2i}\overline{A}%
_{2j}\left(  ||p_{i}||^{2}+||p_{j}||^{2}\right)  \\
  \quad+c\sum_{i<j}\overline{A}_{3i}\overline{A}_{3j}\left(  ||p_{i}%
||^{2}+||p_{j}||^{2}\right)  \\
  =\left(  C\gamma\right)  _{1}\left\Vert p_{1}\right\Vert ^{2}+(C\gamma
)_{2}\left\Vert p_{2}\right\Vert ^{2}+\left(  C\gamma\right)  _{3}\left\Vert
p_{2}\right\Vert ^{2}.
\end{eqnarray*}

\end{proof}

Now we give a lemma which will be the main tool in the construction of
$\gamma_{l}$ from Definition \ref{def:cone-conditions}.

\begin{lemma}
\label{lem:matrix-bounds-cc}Let $U_{i_{1}},U_{i_{2}}\subset\Lambda$ and let
$N$ be a ch-set $N\subset\mathrm{dom}(f_{i_{2}i_{1}}).$ Let $\varepsilon>0$ be
a small number. Let $C$ be defined by (\ref{eq:C-def}) and $\varepsilon>0$.
Assume that $C$ is invertible and define
\begin{equation}
G_{i_{2}i_{1}}=C^{-1}.\label{eq:Gmatrix}%
\end{equation}
If for $\gamma^{\prime}=(a,b,c):=G_{i_{2}i_{1}}\gamma+(\varepsilon
,\varepsilon,\varepsilon),$ we have $a>0,$ and $b,c<0$ then for any
$q_{1},q_{2}\in N$%
\[
Q_{\gamma^{\prime}}(f_{i_{2}i_{1}}(q_{1})-f_{i_{2}i_{1}}(q_{2}))>Q_{\gamma
}(q_{1}-q_{2}).
\]

\end{lemma}

\begin{proof}
For
\[
A:=\int_{0}^{1}Df_{i_{2}i_{1}}(q_{2}+t(q_{1}-q_{2}))dt\in\lbrack
df_{i_{2}i_{1}}(\mathrm{dom}f_{i_{2}i_{1}})]
\]
applying Lemma \ref{lem:matrix-bounds} gives%
\begin{eqnarray*}
Q_{\gamma^{\prime}}(f_{i_{2}i_{1}}(q_{1})-f_{i_{2}i_{1}}(q_{2}))  &
>&Q_{G_{i_{2}i_{1}}\gamma}(f_{i_{2}i_{1}}(q_{1})-f_{i_{2}i_{1}}(q_{2}))\\
&  \geq& Q_{CG_{i_{2}i_{1}}\gamma}(q_{1}-q_{2})\\
&  =&Q_{\gamma}(q_{1}-q_{2}).
\end{eqnarray*}

\end{proof}

\begin{example}
We return to Example \ref{ex:cone-cond-def}. The cones $\gamma_{l}$ follow
from Lemma \ref{lem:matrix-bounds-cc} as $\gamma_{0}=(1,-1)$ and $\gamma
_{l+1}=(1+\varepsilon,(1+\varepsilon)^{-1})\cdot G_{i_{l+1}i_{l}}\gamma_{l}$
with $G_{i_{l+1}i_{l}}=\mathrm{diag}\left(  \frac{1}{\left(  A_{11}%
^{l+1}\right)  ^{2}},\frac{1}{\left(  A_{22}^{l+1}\right)  ^{2}}\right)  ,$
where $\cdot$ stands for the scalar product.
\end{example}

\subsection{Setting up local maps\label{sec:local-maps}}

In this section we shall introduce conditions, which would ensure that the
assumptions from Section \ref{sec:cover-cc-norm-hyp} hold. Below we give a
Lemma which will ensure conditions (\ref{eq:from-g1-to-g0}) and
(\ref{eq:h'-hor-disc}).

Let us note that in some cases conditions (\ref{eq:from-g1-to-g0}) and
(\ref{eq:h'-hor-disc}) will follow from easier arguments or directly from the
setup of the problem. Such is the case in our example from Section
\ref{sec:examples}.

\begin{lemma}
Let $\mathbf{m}>1,$ $\Delta>0$ and $\rho>\sqrt{\frac{\mathbf{a}_{0}%
}{-\mathbf{c}_{0}}}r+\Delta$. Assume that

\begin{enumerate}
\item \label{itm:lem-setup1}for any $\iota\in I$ and any $\lambda\in U_{\iota
}$ there exists a $\boldsymbol{\lambda}_{\kappa}$ such that $(\iota,\kappa)\in
J$ and
\begin{equation}
\left\Vert \eta_{\iota}(\lambda)-\eta_{\iota}(\boldsymbol{\lambda}_{\kappa
})\right\Vert <\Delta,\label{eq:lambda-Delta}%
\end{equation}

\item \label{itm:lem-setup2}for any $\theta\in B_{c}$ and any $i\in I$ there
exists an $\iota\in I$ such that
\begin{eqnarray}
\overline{B}_{c}\left(  \theta,\sqrt{\frac{\mathbf{a}_{1}}{-\mathbf{c}_{1}}%
}r\right)  \cap\overline{B}_{c}   \subset\mathrm{dom}\left(  \eta_{\iota
}\circ\eta_{i}^{-1}\right)  ,\label{eq:dom-cond}\\
\eta_{\iota}\circ\eta_{i}^{-1}(\theta)   \in B_{c}\left(  0,R-\rho
-\Delta\right)  .\label{eq:image-cond}%
\end{eqnarray}
For $C_{\iota\,i}$ defined as in (\ref{eq:C-def}), constructed for
$[d(\tilde{\eta}_{\iota}\circ\tilde{\eta}_{i}^{-1})(\mathrm{dom}(\eta_{\iota
}\circ\eta_{i}^{-1}))]\ $we assume that it is invertible and also that for
$\gamma=(a,b,c)=C_{\iota\,i}^{-1}\boldsymbol{\gamma}_{1}$ we have%
\begin{equation}
\mathbf{a}_{0}>\mathbf{m}a,\quad\mathbf{b}_{0}>\mathbf{m}b,\quad\mathbf{c}%
_{0}>\mathbf{m}c.\label{eq:to-gamma0-bounds}%
\end{equation}

\end{enumerate}

If assumptions \ref{itm:lem-setup1}, \ref{itm:lem-setup2} hold, then for any
horizontal disc $\mathbf{h}$ in a ch-set with cones $(M,\boldsymbol{\gamma
}_{1})$ and for any $i\in I$ there exists $\left(  \iota,\kappa\right)  \in J$
such that $\mathbf{h}(\overline{B}_{u}(0,r))\subset\mathrm{dom}(\tilde{\eta
}_{\iota}\circ\tilde{\eta}_{i}^{-1}).$ Also for any $q_{1},q_{2}$ in
$\mathrm{dom}(\tilde{\eta}_{\iota}\circ\tilde{\eta}_{i}^{-1})$ such that
$Q_{\boldsymbol{\gamma}_{1}}(q_{1}-q_{2})>0$ we have (\ref{eq:from-g1-to-g0}).
Furthermore condition (\ref{eq:h'-hor-disc}) holds.
\end{lemma}

\begin{proof}
Let $\mathbf{h}$ be a horizontal disc in a ch-set with cones
$(M,\boldsymbol{\gamma}_{1}).$ Take $\theta_{0}=\pi_{\theta}(\mathbf{h}(0)).$
For any $x\in\overline{B}_{u}(0,r)$ we have $Q_{\gamma_{1}}(\mathbf{h}%
(x)-\mathbf{h}(0))\geq0,$ which implies that%
\[
\mathbf{a}_{1}r^{2}\geq\mathbf{a}_{1}\left\Vert \pi_{x}(\mathbf{h}%
(x)-\mathbf{h}(0))\right\Vert ^{2}\geq-\mathbf{c}_{1}\left\Vert \pi_{\theta
}(\mathbf{h}(x))-\theta_{0}\right\Vert ^{2},
\]
hence $\pi_{\theta}(\mathbf{h}(\overline{B}_{u}(0,r)))\subset\overline{B}%
_{c}(\theta_{0},\sqrt{\frac{\mathbf{a}_{1}}{-\mathbf{c}_{1}}}r)\cap
\overline{B}_{c}.$ Taking $\iota$ from assumption 2. for $\theta=\theta_{0},$
condition (\ref{eq:dom-cond}) implies that $\mathbf{h}(\overline{B}%
_{u}(0,r))\subset\mathrm{dom}(\tilde{\eta}_{\iota}\circ\tilde{\eta}_{i}^{-1})$
and also
\begin{equation}
\left\Vert \eta_{\iota}\circ\eta_{i}^{-1}(\theta_{0})\right\Vert
<R-\rho-\Delta.\label{eq:theta0-bound}%
\end{equation}
Take now any $q_{1},q_{2}$ in $\mathrm{dom}(\tilde{\eta}_{\iota}\circ
\tilde{\eta}_{i}^{-1})$ such that $Q_{\boldsymbol{\gamma}_{1}}(q_{1}%
-q_{2})>0.$ Applying (\ref{eq:to-gamma0-bounds}) and Lemma
\ref{lem:matrix-bounds-cc} gives%
\begin{eqnarray}
Q_{\boldsymbol{\gamma}_{0}}(\tilde{\eta}_{\iota}\circ\tilde{\eta}_{i}%
^{-1}(q_{1})-\tilde{\eta}_{\iota}\circ\tilde{\eta}_{i}^{-1}(q_{2})) &
>&\mathbf{m}Q_{\gamma}(\tilde{\eta}_{\iota}\circ\tilde{\eta}_{i}^{-1}%
(q_{1})-\tilde{\eta}_{\iota}\circ\tilde{\eta}_{i}^{-1}(q_{2}))\nonumber\\
&  \geq&\mathbf{m}Q_{C_{\iota\,i}C_{\iota\,i}^{-1}\boldsymbol{\gamma}_{1}%
}(q_{1}-q_{2})\nonumber\\
&  =&\mathbf{m}Q_{\boldsymbol{\gamma}_{1}}(q_{1}-q_{2})\label{eq:temp-cc-1}\\
&  >&0,\nonumber
\end{eqnarray}
which proves (\ref{eq:from-g1-to-g0}). Applying the bound in
(\ref{eq:temp-cc-1}) for $q_{1}=\mathbf{h}(x_{1}),$ $q_{2}=\mathbf{h}(x_{2})$
gives
\begin{equation}
Q_{\boldsymbol{\gamma}_{0}}(\mathbf{h}^{\prime}(x_{1})-\mathbf{h}^{\prime
}(x_{2}))\geq0,\label{eq:h'-cone-cond}%
\end{equation}
which means that to prove (\ref{eq:h'-hor-disc}) it is sufficient to show that
$\mathbf{h}^{\prime}(\overline{B}_{u}(0,r))\subset M_{\iota,\kappa}$ for some
$\kappa.$ Let $\lambda=\eta_{i}^{-1}(\theta_{0}).$ We now take $\kappa$ from
assumption 1. For any $x\in\overline{B}_{u}(0,r),$ by (\ref{eq:h'-cone-cond})
we have%
\begin{eqnarray*}
\mathbf{a}_{0}r^{2} &  \geq &\mathbf{a}_{0}\left\Vert \pi_{x}(\mathbf{h}%
^{\prime}(x)-\mathbf{h}^{\prime}(0))\right\Vert ^{2}\\
&  \geq&-\mathbf{c}_{0}\left\Vert \pi_{\theta}(\mathbf{h}^{\prime
}(x)-\mathbf{h}^{\prime}(0))\right\Vert ^{2}\\
&  =&-\mathbf{c}_{0}\left\Vert \pi_{\theta}(\mathbf{h}^{\prime}(x))-\eta
_{\iota}(\lambda)\right\Vert ^{2}.
\end{eqnarray*}
This means that
\[
\pi_{\theta}(\mathbf{h}^{\prime}(\overline{B}_{u}(0,r)))\subset\overline
{B}_{c}(\eta_{\iota}(\lambda),r\sqrt{\frac{\mathbf{a}_{0}}{-\mathbf{c}_{0}}})
\]
hence%
\begin{eqnarray*}
\left\Vert \pi_{\theta}(\mathbf{h}^{\prime}(x))-\eta_{\iota}%
(\boldsymbol{\lambda}_{\kappa})\right\Vert  &  \leq &\left\Vert \pi_{\theta
}(\mathbf{h}^{\prime}(x))-\eta_{\iota}(\lambda)\right\Vert +\left\Vert
\eta_{\iota}(\lambda)-\eta_{\iota}(\boldsymbol{\lambda}_{\kappa})\right\Vert
\\
&  <&r\sqrt{\frac{\mathbf{a}_{0}}{-\mathbf{c}_{0}}}+\Delta\\
&  <&\rho,
\end{eqnarray*}
which gives $\mathbf{h}^{\prime}(\overline{B}_{u}(0,r))\subset M_{\iota
,\kappa}.$ What needs to be verified last is whether $M_{\iota,\kappa}%
\subset\mathbf{B.}$ From our construction $\pi_{\theta}M_{\iota,\kappa
}=\overline{B}_{c}(\eta_{\iota}(\boldsymbol{\lambda}_{\kappa}),\rho).$ For
$\theta\in\overline{B}_{c}(\eta_{\iota}(\boldsymbol{\lambda}_{\kappa}),\rho),$
using (\ref{eq:lambda-Delta}) and (\ref{eq:theta0-bound})%
\begin{eqnarray*}
\left\Vert \theta\right\Vert  &  \leq &\left\Vert \theta-\eta_{\iota
}(\boldsymbol{\lambda}_{\kappa})\right\Vert +\left\Vert \eta_{\iota
}(\boldsymbol{\lambda}_{\kappa})-\eta_{\iota}(\lambda)\right\Vert +\left\Vert
\eta_{\iota}(\lambda)\right\Vert \\
&  =&\left\Vert \theta-\eta_{\iota}(\boldsymbol{\lambda}_{\kappa})\right\Vert
+\left\Vert \eta_{\iota}(\boldsymbol{\lambda}_{\kappa})-\eta_{\iota}%
(\lambda)\right\Vert +\left\Vert \eta_{\iota}\circ\eta_{i}^{-1}(\theta
_{0})\right\Vert \\
&  <&\rho+\Delta+(R-\rho-\Delta),
\end{eqnarray*}
hence $\pi_{\theta}M_{\iota,\kappa}\subset B_{c}.$
\end{proof}

\subsection{Normally hyperbolic manifolds from bounds on derivatives}

In Section \ref{sec:ver-cover} we have shown how covering relations from the
chain (\ref{eq:cover-sequence}) can be constructed using bounds on derivatives
of local maps. In Section \ref{sec:ver-cc} we have shown how the cones can be
set up, using bounds on derivatives of local maps, so that the condition
(\ref{eq:cc-local-iterates}) holds. Here we shall combine these results
together in Theorem \ref{th:main}.

We shall use the notations $T_{i_{2}i_{1}}$ and $G_{i_{2}i_{1}}$ introduced in
Sections \ref{sec:ver-cover}, \ref{sec:ver-cc} through equations
(\ref{eq:Tmatrix-bounds}), (\ref{eq:C-def}) and (\ref{eq:Gmatrix}). We will
also assume that the assumptions from Section \ref{sec:cover-cc-norm-hyp}
hold. Here we introduce a definition which contains conditions which can be
verified using computer assistance. We will later show that the conditions
imply cone conditions.

\begin{definition}
\label{def:Forward-bounds}Assume that for any $(\iota_{0},\kappa_{0})\in J$
there exists an $n\in\mathbb{N,}$ a sequence $\iota_{0}=i_{0},i_{1}%
,\ldots,i_{n}=\iota_{1}$ and $\kappa_{1}$ such that $(\iota_{1},\kappa_{1})\in
J$ and for
\begin{eqnarray*}
q^{m}  &  = &(x^{m},y^{m},\theta^{m}):=f_{i_{m}i_{m-1}}\circ\ldots\circ
f_{i_{1}i_{0}}(0,0,\eta_{i_{0}}(\boldsymbol{\lambda}_{\kappa_{0}})),\\
R^{m}  &  = &(r_{u}^{m},r_{s}^{m},r_{c}^{m}):=T_{i_{m}i_{m-1}}\circ\ldots\circ
T_{i_{1}i_{0}}(r,r,\rho),\\
\gamma^{m}  &  := &(a^{m},b^{m},c^{m}):=G_{i_{m}i_{m-1}}\circ\ldots\circ
G_{i_{1}i_{0}}\gamma_{0}%
\end{eqnarray*}
with $m\leq n$ we have%
\begin{eqnarray}
r_{u}^{m}+\left\Vert x^{m}\right\Vert    < 1,\quad\quad r_{s}^{m}+\left\Vert
y^{m}\right\Vert <1,\quad\quad r_{c}^{m}+\left\Vert \theta^{m}\right\Vert
<1,\nonumber\\
r_{u}^{n}   > r+\left\Vert x^{n}\right\Vert ,\quad\quad r_{s}^{n}+\left\Vert
y^{n}\right\Vert <r, \label{eq:radius-bounds}%
\end{eqnarray}
and
\[%
\begin{array}
[c]{lll}%
a^{m}>0, & 0>b^{m}, & 0>c^{m},\\
a^{n}>\mathbf{a}_{1},\quad & b^{n}>\mathbf{b}_{1},\quad & c^{n}>\mathbf{c}%
_{1}.
\end{array}
\]
Then we say that $f$ \emph{satisfies forward bounds.}
\end{definition}

\begin{remark}
To verify that $f$ satisfies forward bounds on needs to compute $q^{m},$
$R^{m}$ and $\gamma^{m}.$ Let us note that in the case of $q^{m}$ it is enough
to obtain bounds on a finite number of successive iterates of a single point.
We therefore do not need to obtain bounds on images of large sets, which in
practise would accumulate large errors. The $R^{m}$ and $\gamma^{m}$ are
constructed using local bounds on derivatives and are easily computable with
computer assistance. Let us also note that to verify forward bounds we do not
need to compute the composition function $f^{n}$ or its derivative (this would
most likely cause big difficulties for high $n$ due to complexity of such
computations and also due to the fact that errors would accumulate quickly).
\end{remark}

\begin{lemma}
\label{lem:cc-from-forw-bounds}If $f$ satisfies forward bounds then $f$
satisfies cone conditions.
\end{lemma}

\begin{proof}
We take any $(\iota_{0},\kappa_{0})\in J$, a sequence $\iota_{0}=i_{0}%
,i_{1},\ldots,i_{n}=\iota_{1}$ and an index $\kappa_{1}$ such that $(\iota
_{1},\kappa_{1})\in J$ from Definition \ref{def:Forward-bounds}. We define
$R^{0}=R_{\varepsilon}^{0}:=(r,r,\rho)$ and%
\begin{eqnarray*}
R_{\varepsilon}^{m}  &  := &T_{i_{m}i_{m-1}}R_{\varepsilon}^{m-1}+(-\varepsilon
,\varepsilon,\varepsilon)\\
N_{m}  &  :=&N(q^{m},R_{\varepsilon}^{m}).
\end{eqnarray*}
By (\ref{eq:radius-bounds}), taking $\varepsilon$ sufficiently small, we will
ensure that $N_{m}\subset\mathbf{B.}$ By Lemma
\ref{lem:matrix-bounds-covering} we obtain $N_{m-1}\overset{f_{i_{m}i_{m-1}}%
}{\Longrightarrow}N_{m}$ for $m=1,\ldots,n$ and $N_{n}\overset{\mathrm{id}%
}{\Longrightarrow}M.$

Now we define $\gamma^{0}=\gamma_{\varepsilon}^{0}:=\boldsymbol{\gamma}_{0}$
and
\[
\gamma_{\varepsilon}^{m}:=G_{i_{m}i_{m-1}}\gamma_{\varepsilon}^{m-1}%
+(\varepsilon,\varepsilon,\varepsilon).
\]
Taking $\varepsilon>0$ small enough and applying Lemma
\ref{lem:matrix-bounds-cc} we obtain (\ref{eq:cc-local-iterates}).
\end{proof}

From now on let us assume that $f$ satisfies forward bounds with
$\boldsymbol{\gamma}_{0}=\boldsymbol{\gamma}_{0}^{\mathrm{forw}}.$

\begin{definition}
Let $\boldsymbol{\gamma}_{0}^{\mathrm{back}}=\left(  \mathbf{a}_{0}%
^{\mathrm{b}},\mathbf{b}_{0}^{\mathrm{b}},\mathbf{c}_{0}^{\mathrm{b}}\right)
\in\mathbb{R}^{3}$ be such that $\mathbf{a}_{0}^{\mathrm{b}},\mathbf{c}%
_{0}^{\mathrm{b}}<0$ and $\mathbf{b}_{0}^{\mathrm{b}}>0.$ We say that $f$
\emph{satisfies backward bounds} if $f^{-1}$ satisfies forward bounds, with
reversed roles of the $x$ and $y$ coordinates.
\end{definition}

\begin{theorem}
Assume that $f$ satisfies forward bounds for $\boldsymbol{\gamma}%
_{0}^{\mathrm{forw}}=\left(  \mathbf{a}_{0}^{\mathrm{f}},\mathbf{b}%
_{0}^{\mathrm{f}},\mathbf{c}_{0}^{\mathrm{f}}\right)  $ and backward bounds
for $\boldsymbol{\gamma}_{0}^{\mathrm{back}}=\left(  \mathbf{a}_{0}%
^{\mathrm{b}},\mathbf{b}_{0}^{\mathrm{b}},\mathbf{c}_{0}^{\mathrm{b}}\right)
.$ If in addition inequality (\ref{eq:ab-ineq}) holds then there exists a
normally hyperbolic invariant manifold in $\mathcal{U},$ together with its
stable and unstable manifolds $W^{s},$ $W^{u}.$
\end{theorem}

\begin{proof}
This follows directly from Lamma \ref{lem:cc-from-forw-bounds} and Theorem
\ref{th:main}.
\end{proof}

%% file: 05-examples.tex


\section{Example of applications\label{sec:examples}}

\begin{figure}[ptb]
\begin{center}
\includegraphics[width=0.45in]{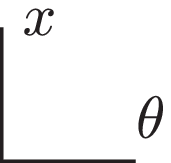}\includegraphics[width=3in]{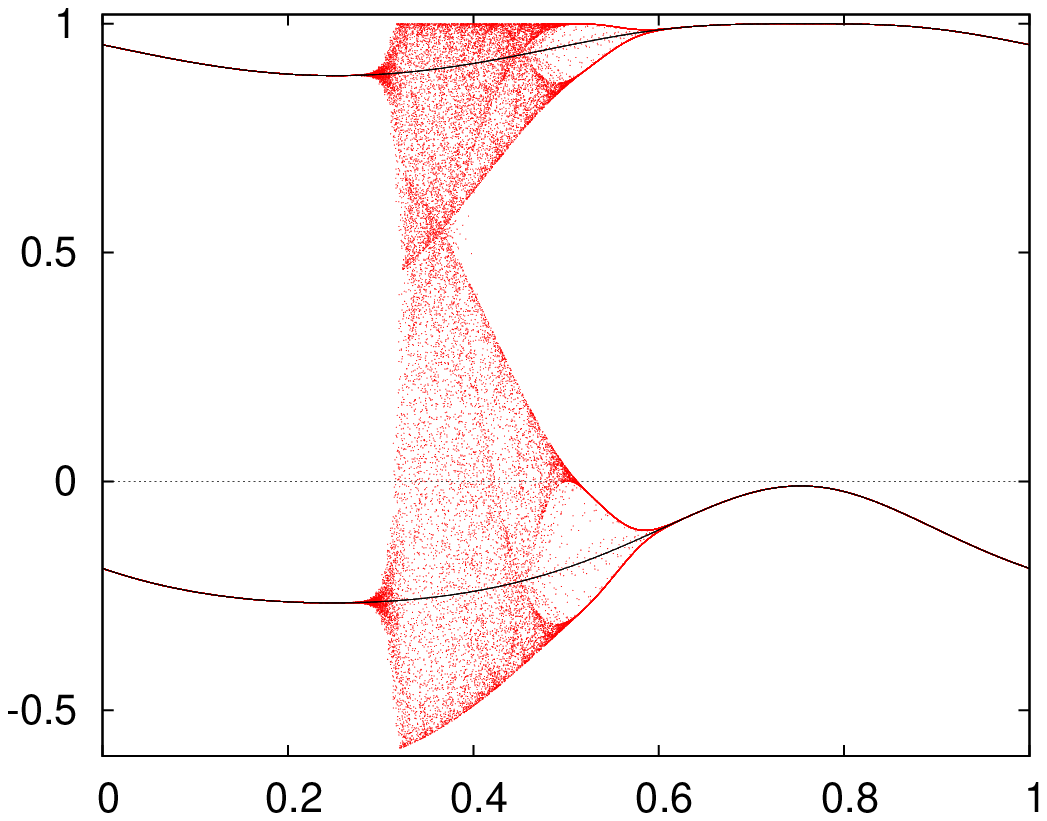}
\end{center}
\caption{Misleading numerical plot of the attractor for $T$, obtained using
double precision (grey), and the true invariant curve computed with 128bit
accuracy (black). }%
\label{fig:misleading200}%
\end{figure}
Consider a driven logistic map
\begin{eqnarray} \label{eq:driven}
& T :S^{1}\times\mathbb{R}\rightarrow S^{1}\times\mathbb{R}, \nonumber \\
& T(\theta,x) =(\theta+\alpha,1-a(\theta)x^2),\quad a(\theta)=
a_0+\eps\sin(2\pi\theta)
\end{eqnarray}
which differs from the well-known logistic map in the fact that the parameter
$a$ has been replaced by $a_0+\eps\sin(2\pi\theta)$ and $\theta$ has a 
quasiperiodic dynamics. Concretely we consider the parameter values 
$a_0=1.31,$ $\eps=0.3$ and $\alpha=\frac{g}{200},$ where $g$ is the golden mean
$g=\frac{\sqrt{5}-1}{2}$, hence the dynamics on the base of the skew-product
is slow. Numerical simulations in double precision (say, with mantissa of 52
binary digits) suggest that the map possesses a chaotic global attractor
(see Figure \ref{fig:misleading200}, grey). We will prove that this guess is not
correct. When the same simulations are done with multiple precision, one can
guess that the attractor consists of two invariant curves (see Figure \ref{fig:misleading200}, black). We will use the
method introduced in the previous sections to prove that $T$ possesses a
contracting invariant manifold and, in particular, that the red plot from Figure
\ref{fig:misleading200} do not shows the true dynamics. The same example was
considered for other values of $\alpha$ and in a non-rigorous way in \cite{BSV}
to illustrate that one has to be careful with the arithmetics in simulations.

\subsection{Explaining the observed behavior}

To explain the reasons of the observed behavior it is worth to mention that in
the example the parameter $a$ of the logistic map ranges in $[a_0-\eps,a_0+\eps]
=[1.01,1.61]$. For that range the attractor starts as a recently created (at
$a=1$) period-2 sink, followed by the full period-doubling cascade. Then one
finds from several-pieces strange attractors to a single piece, interrupted by
some periodic sinks and its corresponding cascades. When $a$ moves with $\theta$
one can question which is the ``averaged'' behavior. In particular the 
period-2 orbit is only attracting until $a=5/4.$

\begin{figure}[htbp]
\begin{center}
\epsfig{file=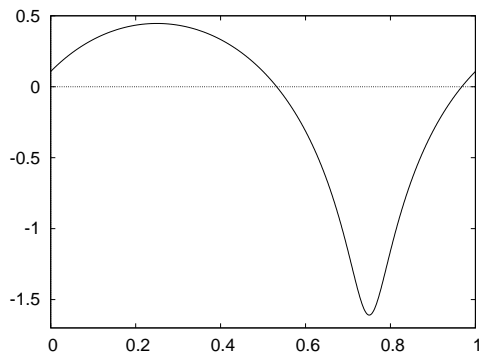,width=7cm}
\end{center}
\caption{The integrand $h(\theta)$ in (\ref{eq:averlyap}) for the parameter
values: $a_0=1.31,\,\eps=0.30.$}
\label{fig:lyapinteg}
\end{figure}

\begin{figure}[htbp]
\begin{center}
\includegraphics[width=0.4in]{figures/coord.eps}\quad
\epsfig{file=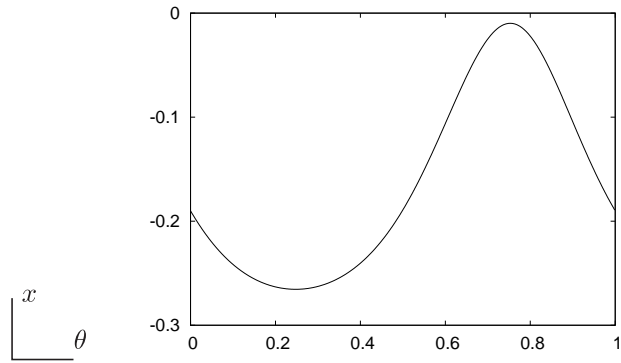,width=7cm}
\end{center}
\caption{The lower part of the attractor, the graph of $x_1(a(\theta))$, for the
parameter values: $a_0=1.31,\,\eps=0.30.$}
\label{fig:attrac.lower}
\end{figure}

To this end we can consider what happens for ``frozen'' values of $a$, denoting
as $T_a$ the corresponding logistic map. The orbit of period two, $x_1(a),
x_2(a)$ is given by the solutions of $x^2-x/a+(1-a)/a^2$. In particular
\be \label{eq:xminus}
x_1(a)=(1-\sqrt{4a-3})/(2a).
\ee
The differential of $T_a^2$ on it is $4(1-a).$ To average with respect to
$\theta$ along the range, and noting that $a-1>0$ for the full range, we have to
consider the average of the Lyapunov exponent given as
\be \label{eq:averlyap}
\frac{1}{2} \int_0^1 \log(4(a_0-1+\eps\sin(2\pi\theta)))d\theta=
\frac{1}{2}\log(2(a_0\!-\!1\!+\!\sqrt{(a_0\!-\!1)^2\!-\!\eps^2}))
\ee
which for $a_0=1.31,\eps=0.3$ gives $\Lambda_\infty\approx -0.12666931.$ The
integrand is shown in Figure \ref{fig:lyapinteg}. For the skew product, assuming
$\alpha\notin\Qset$ and sufficiently small the two curves which form the
attractor, as will be proved later, are very close to the curves $x_1(a),x_2(a)$
of the frozen system. Figure \ref{fig:attrac.lower} displays the lower one. Also
the Lyapunov exponent of the driven map with $\alpha=g/N, N=200$, computed using
$10^5$ iterates after a transient also of $10^5$ iterates is $\Lambda_{200}
\approx -0.12680$. Using other values of $N$, like $100, 400, 800, 1600$ the
respective values $\Lambda_N$ obtained are $-0.12725,\,-0.12670,\,-0.126696,\,
-0.126689$, tending to the limit $\Lambda_\infty$.

The numerical difficulties are easy to understand. To compute the Lyapunov
exponents, starting at a point $x_0$ and an initial vector $v_0=1$ and setting
$S_0=0$ we compute recurrently 
\[ \hat{v}_{j+1}=DT_a(x_j)(v_j),\quad x_{j+1}=T_a(x_j),\quad 
n_{j+1}=|\hat{v}_{j+1}|, \quad v_{j+1}=\hat{v}_{j+1}/n_{j+1},\]
\[ S_{j+1}=S_j+\log(n_{j+1}).  \]
The values $S_j$ are denoted as Lyapunov sums and the average slope as a 
function of $j$ (if it exists) gives the Lyapunov exponent $\Lambda.$ For
details and generalisations see, e.g., \cite{S01} and \cite{LSSW} and references
therein.

\begin{figure}[htbp]
\begin{center}
\begin{tabular}{cc}
\epsfig{file=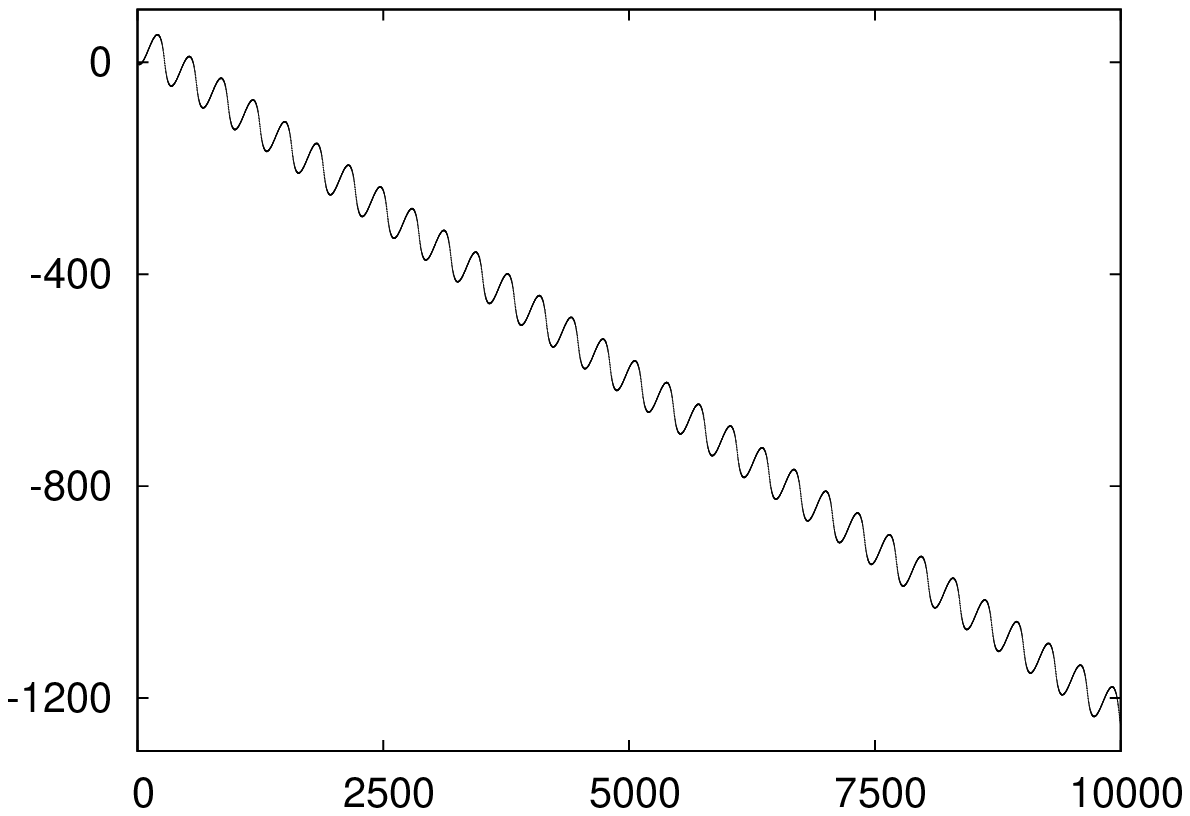,width=5.5cm} &
\epsfig{file=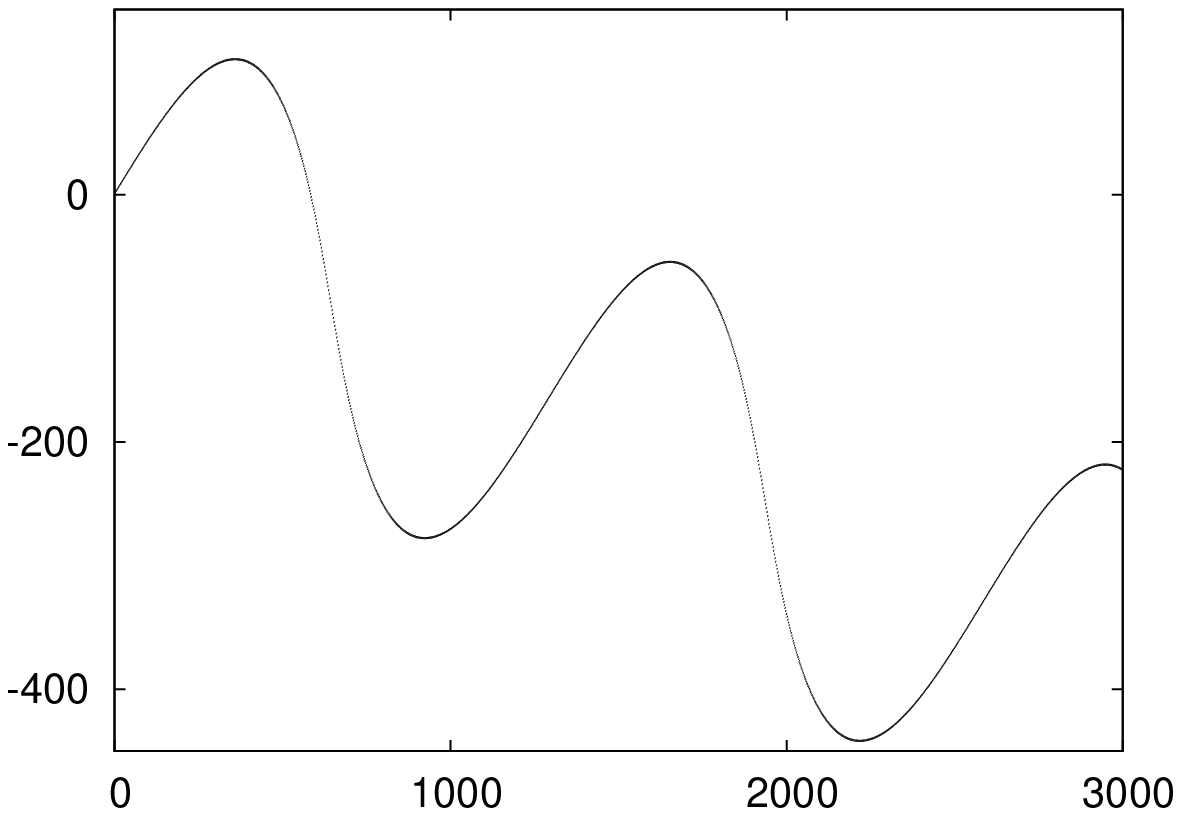,width=5.5cm}
\end{tabular}
\end{center}
\caption{Oscillations of the Lyapunov sums. Left: the Lyapunov sums for $N=200$.
Right: some initial sums for $N=1600.$ Parameter values: $a_0=1.31,\,\eps=0.30$
and $\alpha=g/N.$}
\label{fig:oscillation}
\end{figure}

Even when $\Lambda$ is negative it can happen that partial sums have strong
oscillations. Given the values of $S_j,j=0,\ldots,k$ let $(S_k)_{\text{min}}$
be the minimum of these values and introduce $O_k=S_k-(S_k)_{\text{min}}$. We define the
maximal oscillation of the Lyapunov sums as $OS=max\{O_k\}$. The Figure 
\ref{fig:oscillation} shows the behavior of $S_j$ for $\alpha=g/200$ and also
some of the initial oscillations for $\alpha=g/1600$. A non-rigorous computation
of $OS$ for $N=100,200,400,800,1600$ with $10^5$ iterates after a transient
gives the values $28.845,\,56.761,\,112.632,\,224.379,\,447.874,$ respectively.
This implies a loss in the number of decimal digits equal to these values
divided by $\log(10)$. In particular, between 24 and 25 digits for $N=200$,
which explains the failure seen in Figure \ref{fig:misleading200}. For small
$\alpha$ the maximal oscillation tends to be
\be \label{eq:oscill}
\frac{1}{\alpha}\int_{\theta_2-1}^{\theta_1} h(\theta)d\theta,
\ee
where $h(\theta)$ is the function which appears as integrand in 
(\ref{eq:averlyap}) and it is extended by periodicity outside $[0,1]$ while
$\theta_1=\frac{3}{4}-\frac{1}{2\pi}\cos^{-1}(0.2),$
$\theta_2=\frac{3}{4}+\frac{1}{2\pi}\cos^{-1}(0.2)$ are the values at which $h$
becomes equal to zero (see Figure \ref{fig:lyapinteg}). The value of the
maximal oscillation in (\ref{eq:oscill}) is $\approx 0.172660185/\alpha$ for
small $\alpha$, that if $\alpha=g/N$ becomes $\approx 0.27937N$ in good
agreement with the previous results.

Using these ideas one can even predict when we shall observe that the attractor
produced by simulations with not enough digits seems to indicate that it is not
a period-2 curve. Assume that we do computations with $d$ decimal digits and 
that in a plot like the one in Figure \ref{fig:misleading200} one can 
distinguish pixels which are a a distance of $10^{-p}$. In our example 
reasonable values of $d,p$ are 16 and 4. This means that from $\theta_2-1$, when
$h$ becomes positive, till some unknown $\theta_d$ when the ``departure'' of the
iterates from the curve become visible, the factor of amplification of errors is
$10^{d-p}$ or, in logarithmic scale $(d-p)\log(10)$. This requires
\[\frac{1}{\alpha}\int_{\theta_2-1}^{\theta_d} h(\theta)d\theta=(d-p)\log(10).\]
In our example one finds $\theta_d\approx 0.258$ in good agreement with the
observed numerics in Figure \ref{fig:misleading200} . In a similar way one can
predict the ``landing'' value $\theta_l$ at which the points seen as chaotic in 
Figure \ref{fig:misleading200} are close enough to the real invariant curves.
As the distance from the chaotic points to the true attractor is of the order
of 1, the condition is now
\[\frac{1}{\alpha}\int_{\theta_1}^{\theta_l} h(\theta)d\theta=p\log(10).\]
For the example one obtains $\theta_l\approx 0.629,$ again in good agreement
with the observed numerics.

This ``delayed'' observation of the expanding and compressing regimes is
similar, but now due to purely numerical reasons, to the delay of bifurcation
that can be observed in systems depending on a parameter which has slow dynamics
(see \cite{NST} and references therein).

\subsection{Some limit cases}

Now we discuss two limit cases. First one is the case in which $a(\theta)$
covers a wide range. Second one aims at describing the differences between
the union of the curves $x_1(a(\theta))$ and $x_2(a(\theta))$ and the true
attractor for $\alpha$ small enough.

\begin{figure}[htbp]
\begin{center}
\includegraphics[width=0.4in]{figures/coord.eps}\quad
\epsfig{file=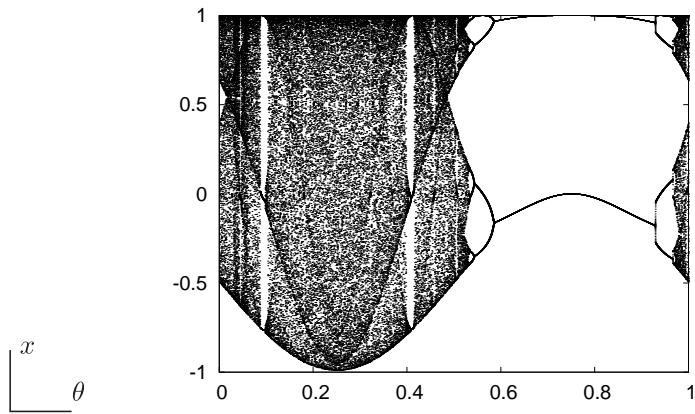,width=8cm}
\end{center}
\caption{Simulations in {\tt double precision} for values of $a_0,\eps$ such
that $a(\theta)$ almost covers the range $(1,2)$ and $\alpha$ very small. See
the text for the numerical values used.}
\label{fig:simbifdia}
\end{figure}

According to (\ref{eq:averlyap}) and assuming that for $\alpha$ sufficiently
small the attractor is close to the union of the curves $x_{1,2}(a(\theta))$
it is enough to take $a_0=1.5-\delta_1,\,\eps=0.5-\delta_1-\delta_2$ with
$0<\delta_2\leq	\delta_1^2$ to have a negative limit averaged Lyapunov exponent
$\Lambda_\infty.$ If $\delta_1$ is small the values of $a$ almost cover the
full range $(1,2)$. The Figure \ref{fig:simbifdia} displays results of the
observed behavior using {\tt double precision} for the values $\delta_1=0.005,\,
\delta_2=10^{-6},\,\alpha=g/60000$.

The figure is reminiscent of the ``bifurcation diagram'' of the logistic map.
In fact, a typical way to compute the diagram consists of taking a sample of
values of $a$, do some transient iterates and display some of the next iterates.
Now the value of $a$ is changed at every step according to (\ref{eq:driven}) 
but very slowly, and the transient is discarded.
From $\theta=3/4$ (for which the minimum of $a(\theta)$ is achieved) to 
$\theta=5/4$ (mod 1) (for which one achieves the maximum) the plot looks like
that diagram, except for the bifurcation delays at the period doublings from 
period 2 to period 4 and successive ones. In the range $\theta\in[1/4,3/4]$ the
reverse situation is seen, but now with much smaller bifurcation delays. The 
authors do not know if, for the present values of the parameter, the attractor
will become close to the union of $x_1(a(\theta))$ and $x_2(a(\theta))$ for
computations done with a huge number of digits.

To look for the expression of the attractor as the union of two smooth curves,
assuming it is of that type, we restrict our attention to the lower part of it,
close to $x_1(a(\theta))$ as given in (\ref{eq:xminus}). In principle it is 
convenient to work with $T^2$ but, as the eigenvalues of $T^2$ along the points
of period are negative, we prefer to work with $T^4$. We look for the attractor
as the graph of a function expanded in powers of $\alpha$
\be \label{eq:graph}
G(\theta)=G_0(\theta)+\alpha G_1(\theta)+\alpha^2G_2(\theta)+\ldots,
\ee
where $G_0(\theta)=x_1(a(\theta))$ is the zeroth order approximation. The map
$T^4(\theta,G(\theta))$ is $\Ocal(\alpha)$ close to the identity. Hence, it can
be approximated by a smooth flow (see \cite{BRS} for proofs, an example of
application  and additional references, as well as \cite{N} for general results)
and the curve we are looking for is a periodic solution of this flow. But we
shall proceed by imposing directly the invariance condition. 

Starting at a point of the form $(\theta,G(\theta))$ and doing four iterations
using the values $a(\theta),a(\theta+\alpha),\,a(\theta+2\alpha),\,
a(\theta+3\alpha)$ we should have
\be \label{eq:invariance}
     T^4(\theta,G(\theta))-(\theta+4\alpha,G(\theta+4\alpha)) =0.
\ee
Given values of $a_0,\,\eps$ it is a cumbersome but elementary task to obtain in
a recurrent way the expressions of $G_1,\,G_2,\ldots$ from 
(\ref{eq:invariance}). It is essential to reduce the dependence in $G_0(\theta)$
using the equation satisfied by $x_1(a)$ to decrease the order of the powers of
$G_0$ which appear to just the first one. We note also that in the computation
of all the terms $G_j$ there appears $16a^2-32a+15=(4a-5)(4a-3)$ in the
denominator, which cancels for $a=5/4$, but a careful examination allows to show
that the factor $4a-5$ is also present in the numerator.

In this way one obtains
\be \label{eq:order1}
G_1(\theta)=\frac{3-2a-(8a-9)/\sqrt{4a-3}}{2a^2(4a-3)}2\pi\eps\cos(2\pi\theta),
\ee
where $a$ stands for $a(\theta)$ as introduced in (\ref{eq:driven}).

\begin{figure}[htbp]
\begin{center}
\begin{tabular}{cc}
\epsfig{file=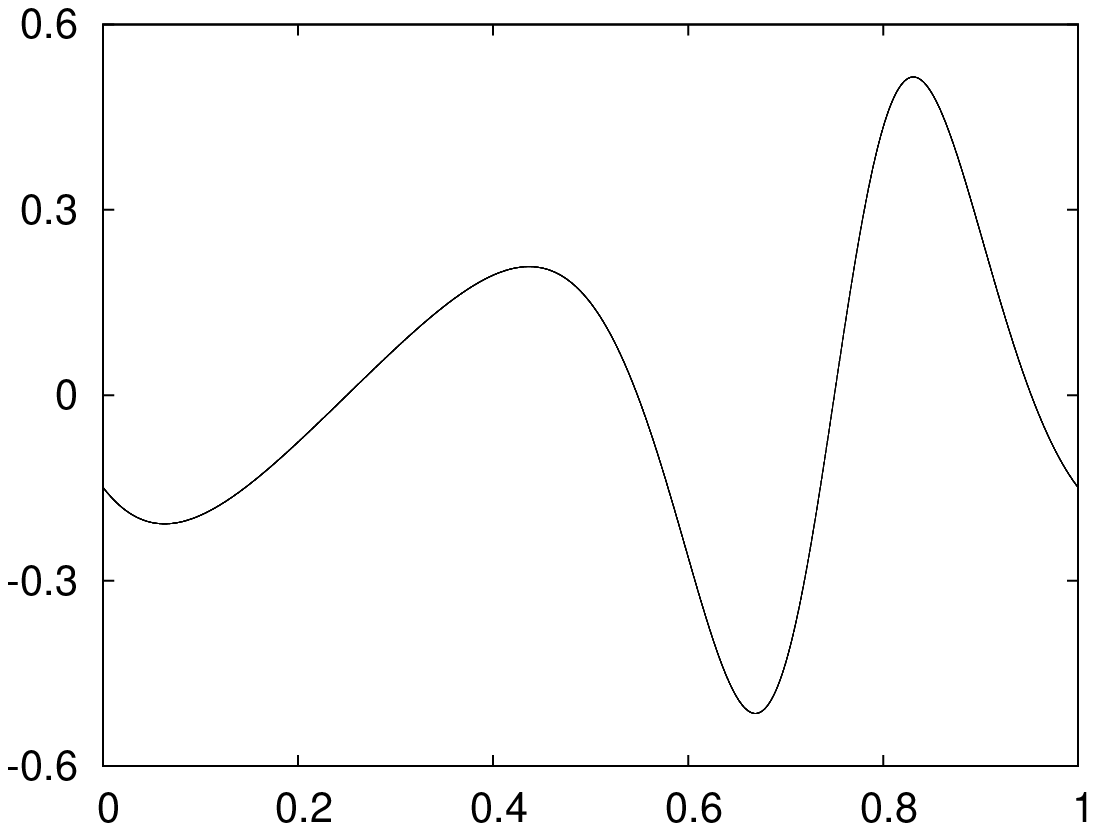,width=5.5cm} &
\epsfig{file=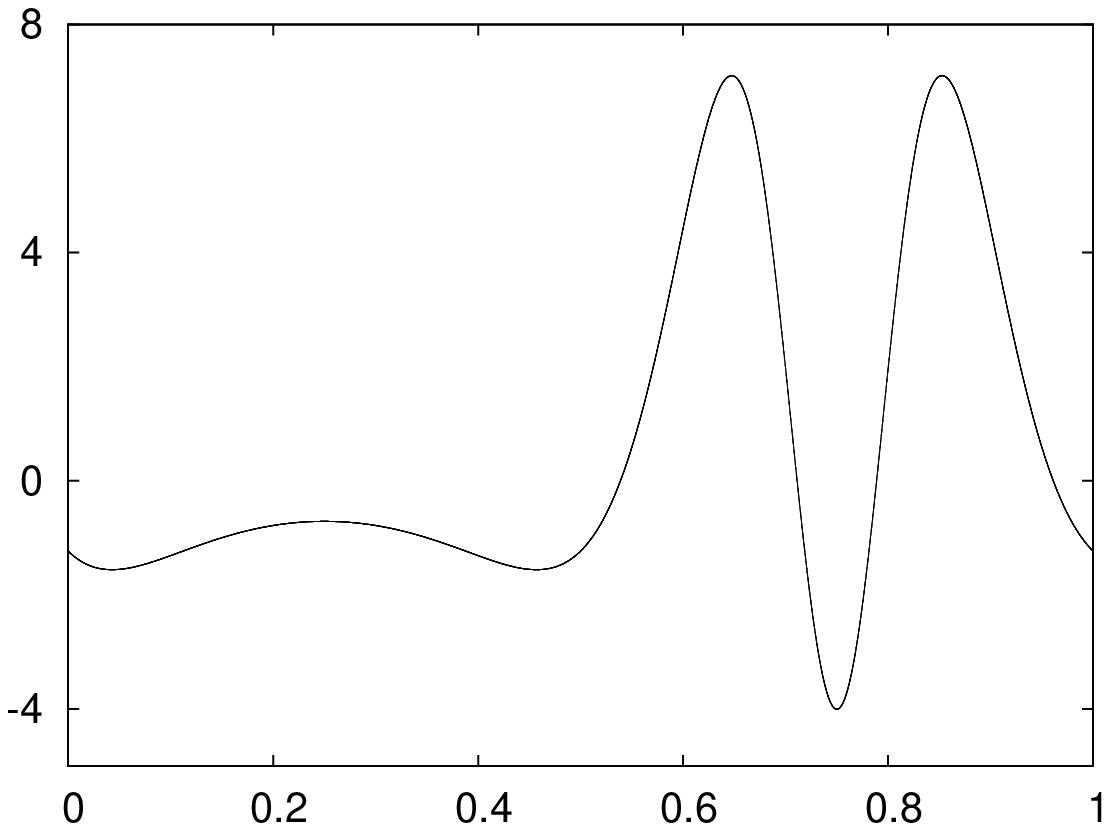,width=5.5cm}
\end{tabular}
\end{center}
\caption{Graphs of $G_1(\theta)$ (left) and $G_2(\theta)$ (right) for
$a_0=1.31,\,\eps=0.30$.}
\label{fig:order1_2}
\end{figure}

The computation of $G_2$ is much more involved. The simplest expression is
given as a rational function depending on $a(\theta),G_0(\theta),G_1(\theta)$
and up to the second derivatives of these functions with respect to $\theta$.
Instead, Figure \ref{fig:order1_2} displays the graph of $G_1$ and $G_2$ for
$a_0=1.31,\,\eps=0.30$. The graph of $G_0(\eps)$ is very close to the attractor
shown in Figure \ref{fig:attrac.lower}.

\begin{figure}[htbp]
\begin{center}
\begin{tabular}{cc}
\epsfig{file=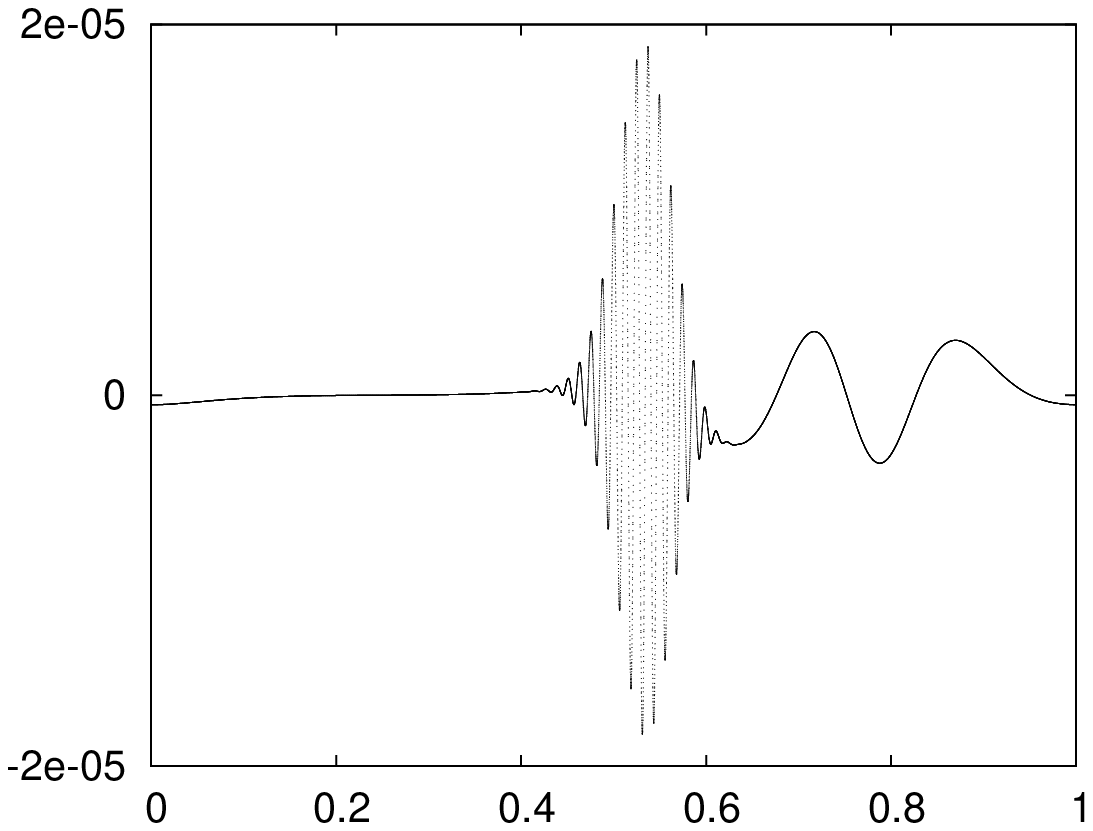,width=5.5cm} &
\epsfig{file=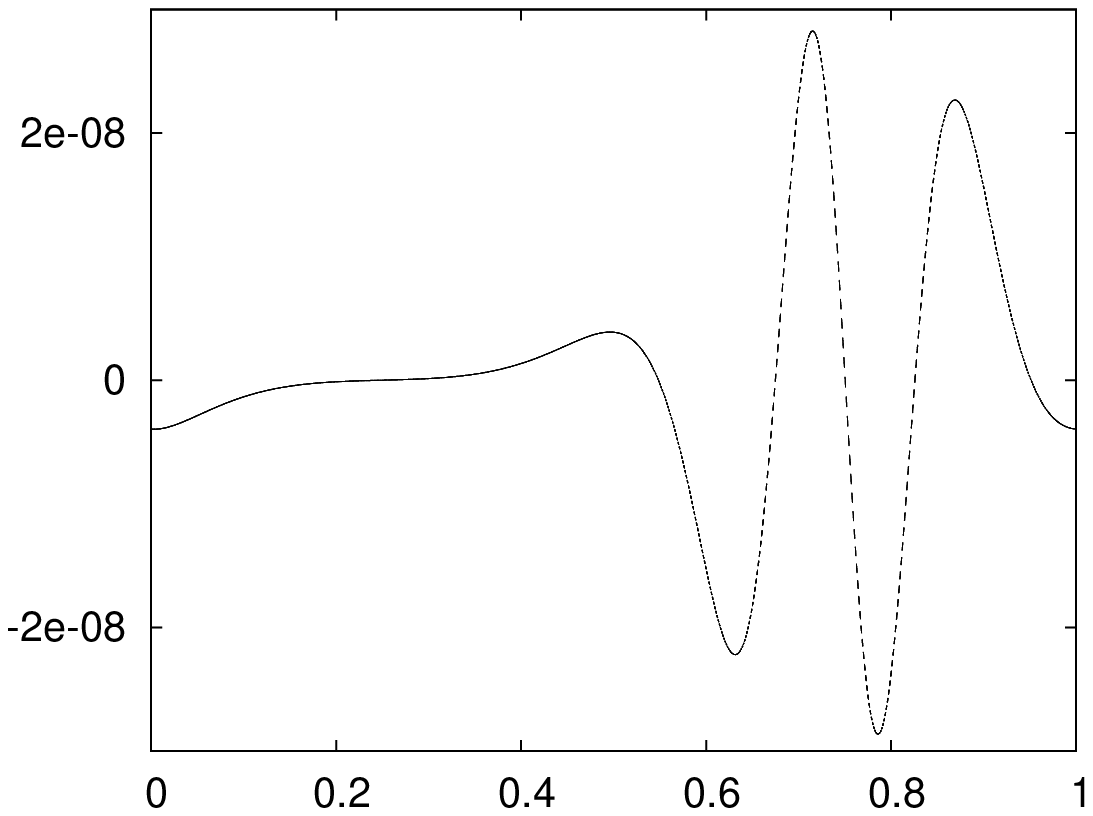,width=5.5cm}
\end{tabular}
\end{center}
\caption{Differences between the attractor and the second order approximation
for $a_0=1.31,\,\eps=0.30$ and $\alpha=g/N$. Left: $N=200$. Right: $N=1000$.}
\label{fig:remainder}
\end{figure}

To see tiny details on the attractor Figure \ref{fig:remainder} displays the
differences between the lower part of the attractor, computed with enough
digits, and the approximation in (\ref{eq:graph}) up to order 2 in $\alpha$.

The left part shows tiny oscillations which were not visible in Figure
\ref{fig:attrac.lower}. They reach a maximum at the value $\theta=\theta_1$ for
which $h(\theta)$ in (\ref{eq:oscill}) changes from positive to negative. As one
can expect the shape of these oscillations is a bump function multiplied by
a periodic function (close to a sinus) with period $4\alpha$. A similar behavior
is observed for many other values of $a_0,\eps$ and $\alpha$. When the 
oscillations start at a larger distance from $\theta_1$ they can amplify is such
a way that the attractor is no longer the union of the two curves. One can
suspect that it becomes a non-chaotic strange attractor (see, e.g., \cite{Greb}
and \cite{Ja}). In contrast, with the same values of $a_0,\eps$ but for $N=1000$
the oscillations are not observed and the very small differences in the plot on 
the right hand side of Figure \ref{fig:remainder} are mainly due to the third
order term in (\ref{eq:graph}).

\subsection{Computer assisted proof of existence of invariant curves}
In this section we apply our method from Sections \ref{sec:geometric}, \ref{sec:covering-ver-ex} to prove that for parameters $a_0=1.31,$ $\eps=0.3$ and $\alpha=\frac{g}{200},$ with $g=\frac{\sqrt{5}-1}{2}$ the map $T$ has an invariant curve. 

\begin{figure}[ptb]
\begin{center}
\includegraphics[width=0.45in]{figures/coord.eps}\hspace{1cm}%
\includegraphics[height = 1.8in]{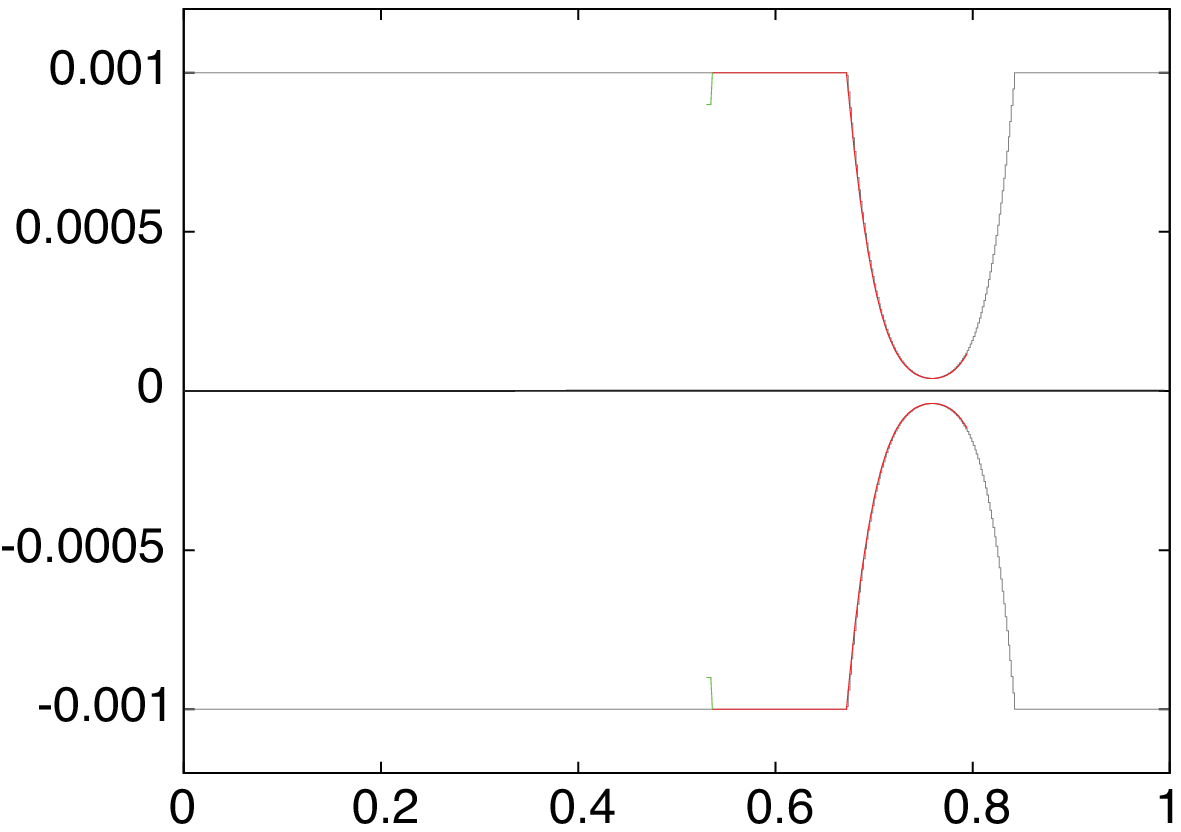}
\end{center}
\caption{Positioning of our ch-sets (green and red) and the set $\mathcal{V}$
(gray) relative to the attractor (on this plot on the $\theta$-axis in black).}%
\label{fig:attr-bounds}%
\end{figure}

Around a neighborhood of the numerical guess for the attractor, the map
$T^{2}$ is locally invertible. This is due to the fact that our curve is
separated from the $x$-axis. For our proof we consider
\[
f=T^{-2}.
\]
The attractor is first computed (nonrigorously) by iterating $T$ forwards in
time. We then choose a set $\mathcal{V}$ arround the attractor (see Figure
\ref{fig:attr-bounds}, gray). For most $\theta$ the set is a $0.001$ radius
neighbourhood of the attractor. Close to the angle $\theta=\frac{3}{4}$ we
choose $\mathcal{V}$ to be tighter, so that we are sure that it lies within
the domain of $f$ (see Figure \ref{fig:attr-bounds}). Our aim is to prove that
inside of $\mathcal{V}$ we have an invariant normally hyperbolic curve of $f.$
The map $f$ is not uniformly expanding in the $x$ direction. Over one part of
the set $\mathcal{V}$ the map $f$ is strongly expanding, elswehere it is
contracting. A part of the expansion region, which we denote as
$\mathcal{U\subset V}$, is depicded in red and green (the green region is on
the left tip of the red region and is poining towards the attractor) on Figure
\ref{fig:attr-bounds}. On this set we place ch-sets $N_{1},\ldots,N_{168}$ of
width $\frac{\alpha}{2}$, starting with $N_{1}$ on the left and finishing with
$N_{168}$ on the right. We shall use a notation%
\[
U_{k,l}=\bigcup_{i=k}^{l}N_{i}.
\]
Our ch-sets are parallelograms. The coordinate $x$ is globally expanding for
$f$ and coordinate $\theta$ is normal (our map does not posses a globally
contracting coordinate $y$). The exits sets $N_{i}^{-}$ for the ch-sets are
the top and bottom edges of the parallelograms. The map $f$ moves the ch-sets
to the left. We distinguish two parts of the set $\mathcal{U}$: the set
$U_{1,4}$ in our plots is denoted in green colour, $U_{5,168}$ is denoted in
red. Since the width of the ch-sets is $\frac{\alpha}{2},$ for $k\in
5,\ldots,168$ we have%
\[
\pi_{\theta}f(N_{k})\subset\pi_{\theta}N_{k-4}.
\]

\begin{figure}[ptb]
\begin{center}
\includegraphics[width=0.4in]{figures/coord.eps}
\includegraphics[height = 1.5in]{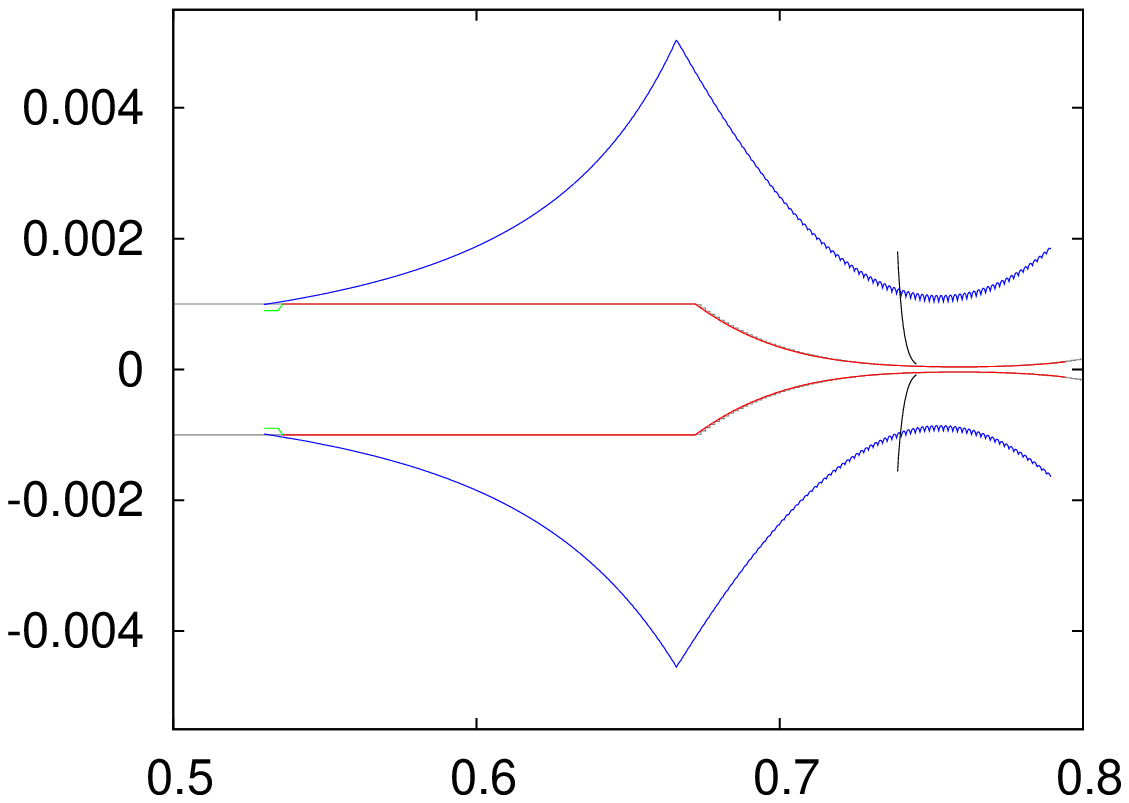}
\includegraphics[height = 1.5in]{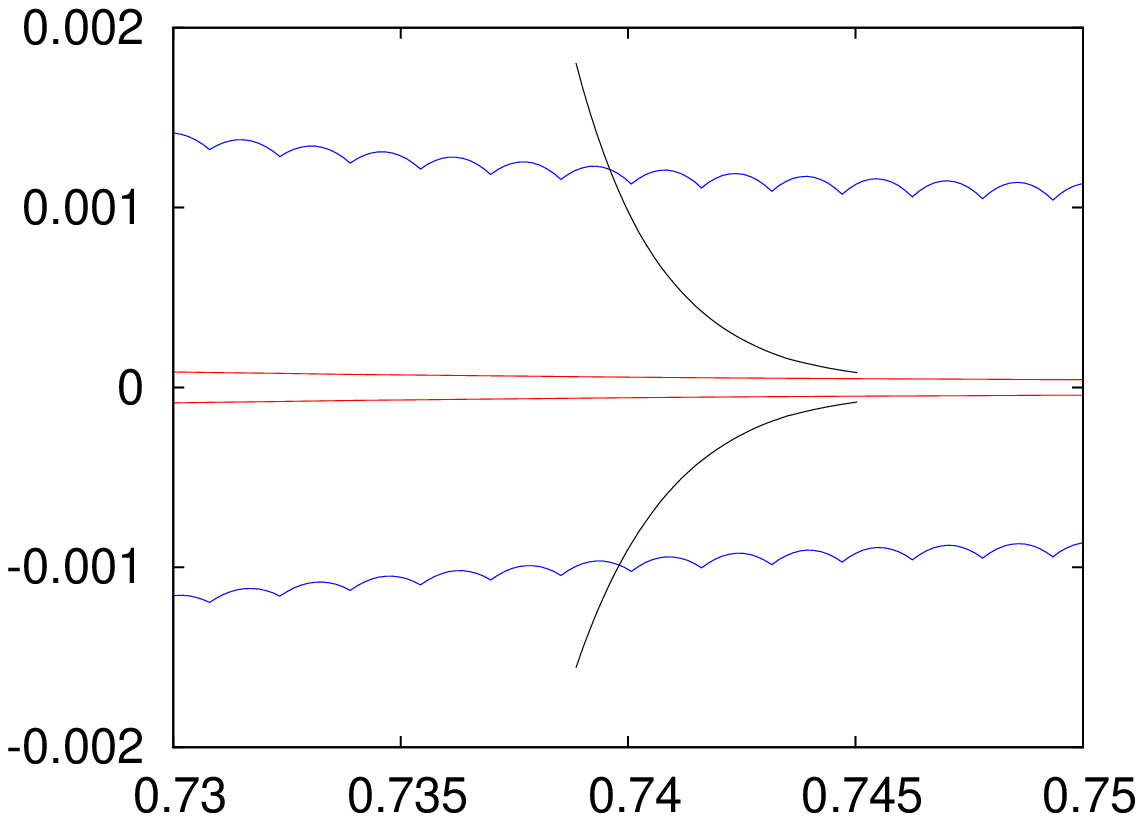}
\end{center}
\caption{The ch-sets: $N_{1}^{-},..N_{4}^{-}$ in green and $N_{5}%
^{-},...,N_{168}^{-}$ in red (plotted relative to the attractor), together
with $f(N_{5}^{-}),...,f(N_{168}^{-})$ in dark blue and $f^{128}(N_{1}%
^{-}),...,f^{128}(N_{4}^{-})$ in black. }%
\label{fig:ch-sets-covering}%
\end{figure}

\begin{figure}[ptb]
\begin{center}
\includegraphics[width=0.4in]{figures/coord.eps}\hspace{1cm}
\includegraphics[ height=1.5in]{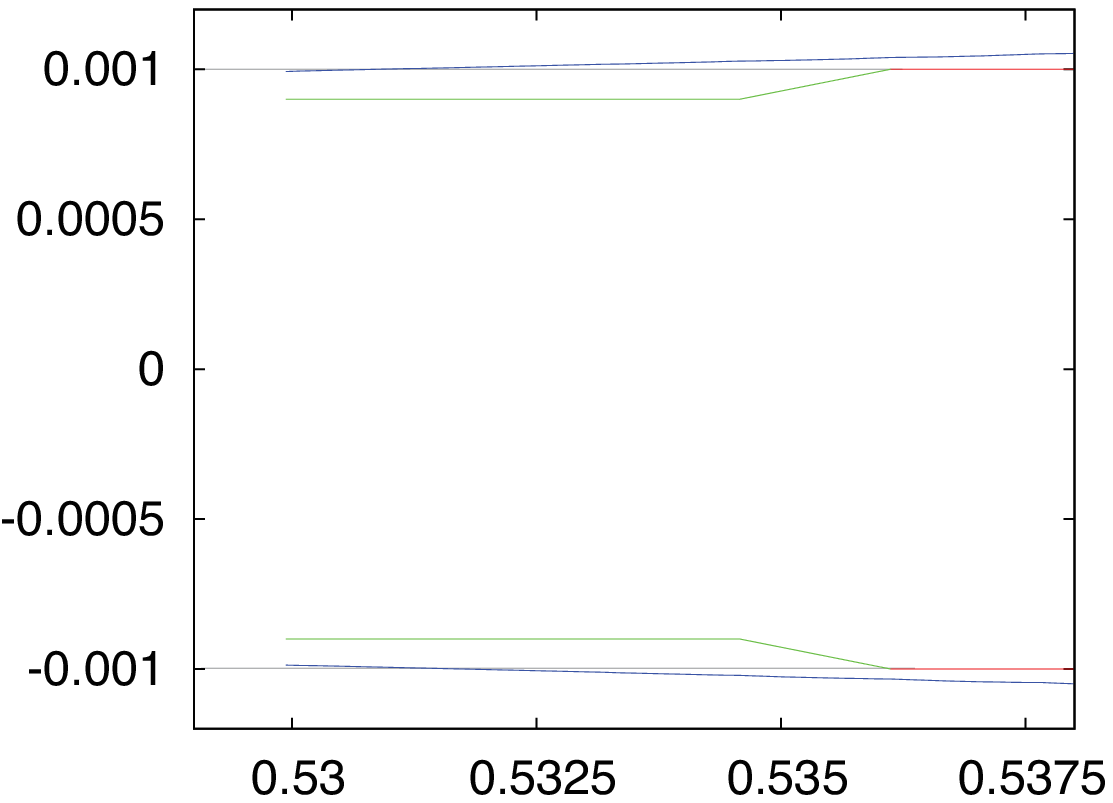}
\end{center}
\caption{Closeup of the covering $N_{i}\overset{f}{\Longrightarrow}N_{i-4}$
for $i=5,6,7,8$ (plotted relative to the attractor).}%
\label{fig:covering2}%
\end{figure}

In Section \ref{sec:covering-ver-ex} we shall show that (see Figures
\ref{fig:ch-sets-covering}, \ref{fig:covering2})%
\begin{equation}
N_{k}\overset{f}{\Longrightarrow}N_{k-4}\quad\text{for }k\in\{5,\ldots,168\},
\label{eq:covering-single}%
\end{equation}
and also that for $i=1,...,4$ (see Figure
\ref{fig:ch-sets-covering})%
\[
N_{i}\overset{f^{128}}{\Longrightarrow}U_{5,168}.
\]

In Section \ref{sec:cc-ver-ex} we show how to verify cone conditions. In
Section \ref{sec:tech-notes} we briefly discribe the tools that were used to
conduct the proof.



\subsubsection{Verification of covering conditions\label{sec:covering-ver-ex}}

To describe how covering conditions are verified we start with a seemingly
unrelated discussion. Consider a polynomial $p:[0,r]\rightarrow\mathbb{R}$ of
degree $n$%
\[
p(\theta)=\sum_{j=0}^{n}a_{j}\theta^{j},
\]
and a function $g:\mathbb{R}\rightarrow\mathbb{R.}$ Using Taylor expansion and
defining two polynomials $\overline{p}$ and $\underline{p}$, of degree $n$%
\begin{eqnarray}
\overline{p}(\theta)  &  = &g\circ p(0)+\sum_{j=1}^{n-1}\left(  \frac{1}%
{j!}\frac{d^{j}\left(  g\circ p\right)  }{d\theta^{j}}(0)\right)  \theta
^{j}\label{eq:p-upper}\\
&&  +\frac{1}{n!}\left(  \frac{d^{n}\left(  g\circ p\right)  }{d\theta^{n}%
}(0)+\frac{1}{n+1}\sup_{v,w\in\lbrack0,r]}\frac{d^{n+1}\left(  g\circ
p\right)  }{d\theta^{n+1}}(v)w\right)  \theta^{n},\nonumber\\
\underline{p}(\theta)  &  =&g\circ p(0)+\sum_{j=1}^{n-1}\left(  \frac{1}%
{j!}\frac{d^{j}\left(  g\circ p\right)  }{d\theta^{j}}(0)\right)  \theta
^{j}\label{eq:p-lower}\\
&&  +\frac{1}{n!}\left(  \frac{d^{n}\left(  g\circ p\right)  }{d\theta^{n}%
}(0)+\frac{1}{n+1}\inf_{v,w\in\lbrack0,r]}\frac{d^{n+1}\left(  g\circ
p\right)  }{d\theta^{n+1}}(v)w\right)  \theta^{n},\nonumber
\end{eqnarray}
for any $\theta\in\lbrack0,r]$ we have%
\begin{equation}
\underline{p}(\theta)\leq g(p(\theta))\leq\overline{p}(\theta).
\label{eq:edge-bounds}%
\end{equation}

For any $i=1,...,168,$ the exit set $N_{i}^{-}$ consists of two lines and can
be expressed using two polynomials (in fact these are affine functions)
$p_{i}^{u},p_{i}^{d}:[0,\frac{\alpha}{2}]\rightarrow\mathbb{R,}$ $p_{i}%
^{d}(\theta)=a_{i,0}^{d}+a_{i,1}^{d}\theta,$ $p_{i}^{u}(\theta)=a_{i,0}%
^{u}+a_{i,1}^{u}\theta$ and a point $q_{i}\in\lbrack0,1),$%
\begin{eqnarray*}
N_{i}^{-}  &  = &N_{d}^{-}\cup N_{u}^{-},\\
N_{i,d}^{-}  &  =&\{(p_{i}^{d}(\theta),q_{i}+\theta)|\theta\in\lbrack
0,\frac{\alpha}{2}]\},\\
N_{i,u}^{-}  &  =&\{(p_{i}^{u}(\theta),q_{i}+\theta)|\theta\in\lbrack
0,\frac{\alpha}{2}]\},
\end{eqnarray*}%
\[
p_{i}^{d}(\theta)<p_{i}^{u}(\theta)\quad\text{for }\theta\in\lbrack
0,\frac{\alpha}{2}].
\]
We will now show how to construct a ch-set $M$ such that%
\begin{equation}
N_{i}\overset{f}{\Longrightarrow}M. \label{eq:ni-covering-1step}%
\end{equation}
We first verify that for any point $(\theta,x)\in N_{i}$ we have
$\frac{\partial f}{\partial x}(x,\theta)<0$. We then take
\begin{equation}
g^{u}(\theta):=f(q_{i}+\theta,p_{i}^{d}(\theta)),\quad g^{d}(\theta
):=f(q_{i}+\theta,p_{i}^{u}(\theta)), \label{eq:g-maps}%
\end{equation}
and construct $p^{u}(\theta)=\overline{p}(\theta)$ using (\ref{eq:p-upper})
and $p^{d}(\theta)=\underline{p}(\theta)$ using (\ref{eq:p-lower}), taking $g$
as functions $g^{u}$ and $g^{d}$ respectively. Formula (\ref{eq:edge-bounds})
guarantees that $f(N_{i,d}^{-})$ lies above the graph of $p^{u}(\theta)$ and
that $f(N_{i,u}^{-})$ lies below the graph of $p^{d}(\theta)$. If we now set%
\begin{eqnarray}
M^{-}  &  =&M_{d}^{-}\cup M_{u}^{-},\nonumber\\
M_{d}^{-}  &  =&\{(p^{d}(\theta),q_{i}-2\alpha+\theta)|\theta\in\lbrack
0,\frac{\alpha}{2}]\},\label{eq:Mi-procedure}\\
M_{u}^{-}  &  =&\{(p^{u}(\theta),q_{i}-2\alpha+\theta)|\theta\in\lbrack
0,\frac{\alpha}{2}]\},\nonumber
\end{eqnarray}
and take $M$ to be the set of points which lie above $M_{d}^{-}$ and below
$M_{u}^{-}$ then (\ref{eq:ni-covering-1step}) holds.

For $i=5,...,168,$ after applying the above procedure to abtain $M$ which is
covered by $N_{i},$ we compute bounds on the images of sets
\begin{equation}
p^{u}\left(  \left[  \frac{j\alpha}{20},\frac{(j+1)\alpha}{20}\right]
\right)  ,\quad p^{d}\left(  \left[  \frac{j\alpha}{20},\frac{(j+1)\alpha}%
{20}\right]  \right)  \quad\text{for }j=0,...,9, \label{eq:boundary-bound}%
\end{equation}
in local coordinates of ch-sets $N_{i-4},$ to verify that we have
(\ref{eq:covering-single}) (subdividing $[0,\frac{\alpha}{2}]$ into ten
intervals turns out to be sufficient for all $i\in\{5,...,168\}$).

For $i=1,...,4$ we need to iterate the procedure (\ref{eq:Mi-procedure}) many
times to obtain a sequence of covering relations%
\[
N_{i}\overset{f}{\Longrightarrow}M_{1}\overset{f}{\Longrightarrow}%
M_{2}\overset{f}{\Longrightarrow}\ldots\overset{f}{\Longrightarrow}%
M_{127}\overset{f}{\Longrightarrow}U_{5,168}.
\]
During our construction we make sure that all sets $M_{k}$ for $k\in
\{1,...,127\}$ lie in $\mathcal{V}$, which readily holds since the sets are
very strongly contracted. Each covering $M_{k}\overset{f}{\Longrightarrow
}M_{k+1}$ holds by construction. Verifying that $M_{127}\overset
{f}{\Longrightarrow}U_{5,168}$ is done analogously to (\ref{eq:boundary-bound}).

In our computer assisted proof we take the degrees of polynomials for the
edges of the sets $M_{k}$ as nine, which means that we need to perform
$C^{10}$ computations. Let us note that computationally this is not as heavy
as might seem, since the $C^{10}$ computations are performed for one
dimensional functions $g^{u}(\theta)$ and $g^{d}(\theta)$ (see
(\ref{eq:g-maps})). The reduction of dimension truly pays off, since the
difference between $C^{10}$ computations in one and two dimensions is substantial.

\begin{figure}[ptb]
\begin{center}
\includegraphics[width=0.4in]{figures/coord.eps}
\includegraphics[width=2in]{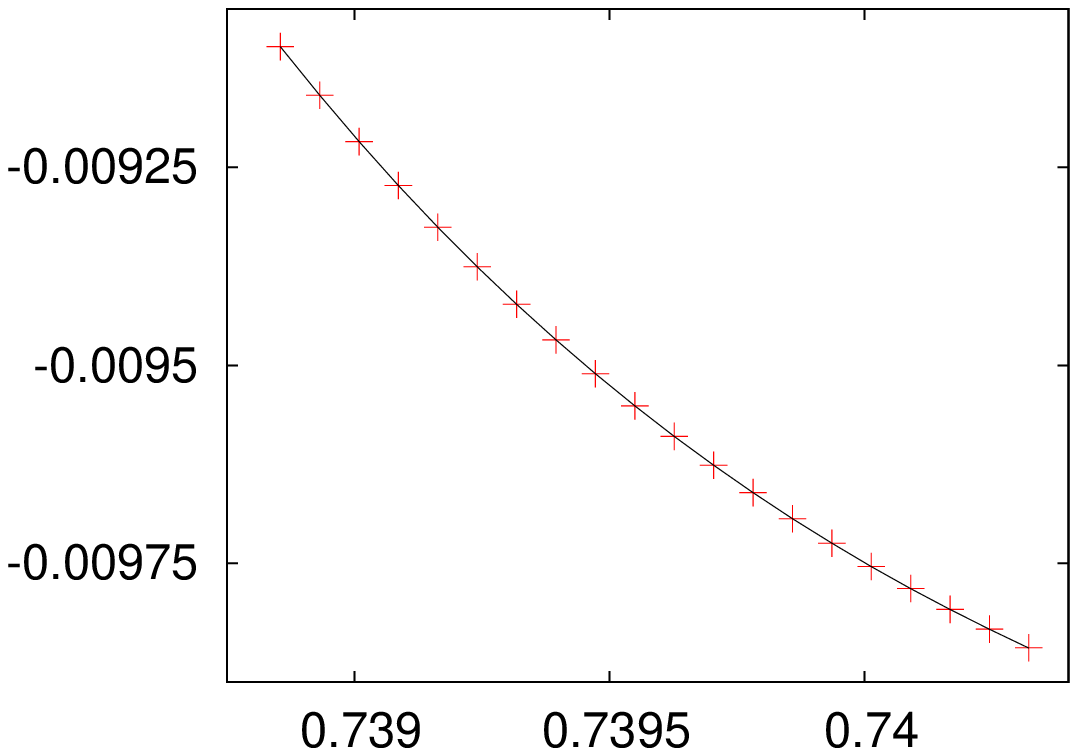}
\includegraphics[width=2in]{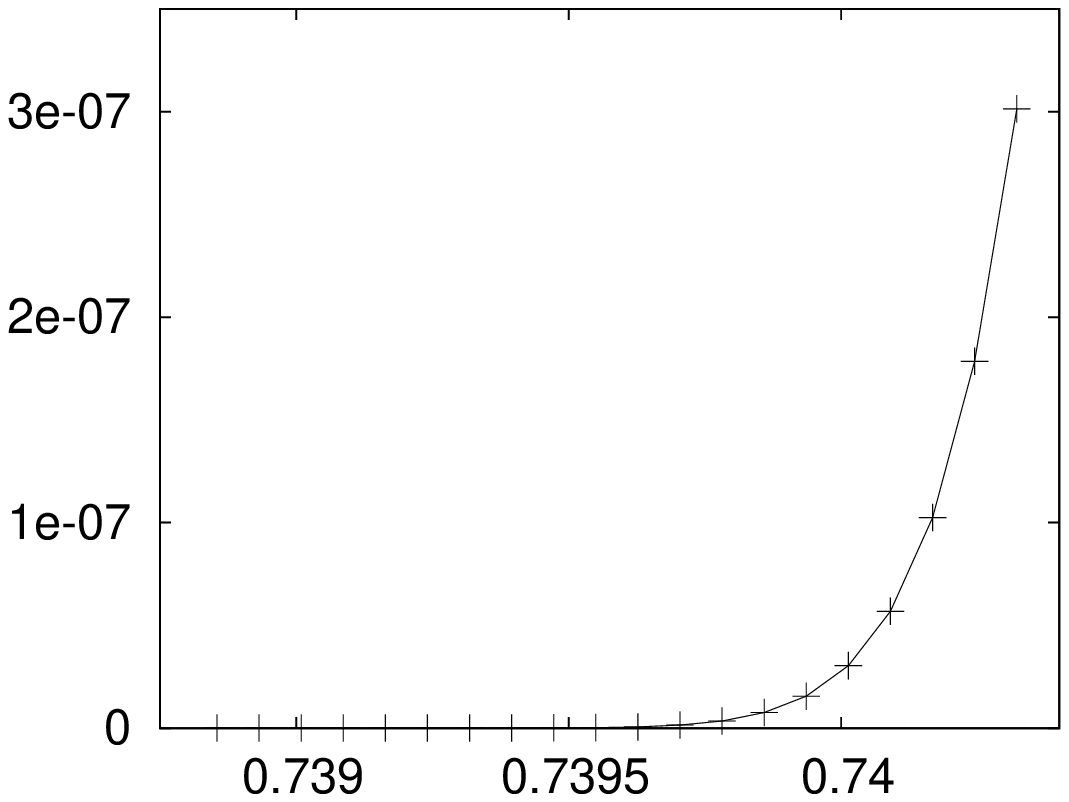}
\end{center}
\caption{Left: Rigorous bound on the image of an edge of one of ch-sets after
128th iterate of the map (in black), together with non-rigorous computations
using multi-precision (red). Right: The difference between rigorous lower
bound and non-rigorous computations. }%
\label{fig:num-difference}%
\end{figure}

The estimates obtained by us are very accurate. In Figure
\ref{fig:num-difference} we give a plot of $M_{128,u}^{-}$, which is the lower
bound estimate of the image of $N_{4,u}^{-}$ after the final step in our
procedure (in black), and compare it with ten points from $N_{4,u}^{-},$
iterated non-rigorously with high precision computations (in red). The curve
lies below the points, as should, but this is impossible to distinguish from
the graph. The right hand side of Figure \ref{fig:num-difference} gives the
plot of the difference of the rigorous lower bound and non-rigorous
computation. They turn out to be very close.

\begin{remark}
The high order computations and multi-precision in current approach seem
essential. The sets $M_{k}$ constructed with our procedure are very strongly
contracted. The distance between the two curves of $M_{k}^{-}$ at the tightest
spot is of order $1.125\times10^{-25}$, which is extremely thin when compared
to the width of the curves $\frac{\alpha}{2}\approx1.545\times10^{-3}$; and
yet, with our $C^{10}$ approach, with little effort we are able to rigorously
keep them apart. Any standard approach, such as performing $C^{0}$
computations on sets or careful linearization with $C^{1}$ techniques through
local coordinates, is likely to fail.
\end{remark}

\begin{remark}
We believe that using a "parallel shooting" type approach it should be
possible to conduct the proof using double precision and $C^{1}$ computations
only (for this we would need an good apriori guess for the position of the
curve). Such approach could produce a rigorous-computer-assisted proof using
double precision of an invariant a curve, which is not detectable numerically with double
computations. This shall be a subject of forthcoming work. 
\end{remark}

\subsubsection{Verification of cone conditions\label{sec:cc-ver-ex}}

To verify cone conditions let us first rescale our coordinates by%
\[
\gamma_{\beta}(\theta,x)=(\beta\theta,x).
\]
Taking $\beta$ sufficiently large, choosing sufficiently many points
$\mathbf{\lambda}_{i}\in\lbrack0,\beta)$ and taking $h_{i}:=\frac{1}{2}\left(
c^{u}\left(  \mathbf{\lambda}_{i}\right)  -c^{d}\left(  \mathbf{\lambda}%
_{i}\right)  \right)  ,$ $q_{i}:=(\mathbf{\lambda}_{i},c^{d}\left(
\mathbf{\lambda}_{i}\right)  +h_{i})$ and $V_{i}:=\mathcal{V}\cap\left(
\lbrack\mathbf{\lambda}_{i}-h_{i},\mathbf{\lambda}_{i}+h_{i}]\times
\mathbb{R}\right)  $ we can construct local maps%
\[
\tilde{\eta}_{i}:V_{i}\rightarrow B_{c}\times B_{u},
\]
for which $\tilde{\eta}_{i}(V_{i}\cap c^{u})=B_{c}\times\{1\},$ $\tilde{\eta
}_{i}(V_{i}\cap c^{d})=B_{c}\times\{-1\}$ and which are arbitrarily close to a
linear map $q\rightarrow\frac{1}{h_{i}}(q-q_{i})$. In these local coordinates,
by taking sufficiently large $\beta$, we have the following bound on
derivatives of local maps (assuming that we choose $i,j$ and $p$ such that
$p\in dom(f_{ij})\neq\emptyset$)%
\begin{eqnarray*}
Df_{ij}  &  =&D\left(  \tilde{\eta}_{i}\circ\gamma_{\beta}\circ f\circ
\gamma_{\beta}^{-1}\circ\tilde{\eta}_{j}^{-1}\right)  (p)\\
&&  \approx\left(
\begin{array}[c]{cc}%
	\frac{1}{h_{i}} & 0\\
	0 & \frac{1}{h_{i}}%
\end{array}
\right)  \left(
\begin{array}[c]{cc}%
	\beta & 0\\
	0 & 1
\end{array}
\right)  \left(
\begin{array}[c]{cc}%
	\frac{df_{1}}{d\theta}(\gamma_{\beta}^{-1}(\tilde{\eta}_{j}^{-1}(p)) & 0\\
	\frac{df_{2}}{d\theta}(\gamma_{\beta}^{-1}(\tilde{\eta}_{j}^{-1}(p)) &
	\frac{df_{2}}{dx}(\gamma_{\beta}^{-1}(\tilde{\eta}_{j}^{-1}(p))
\end{array}
\right) \\
&&  \left(
\begin{array}
[c]{cc}%
\beta^{-1} & 0\\
0 & 1
\end{array}
\right)  \left(
\begin{array}
[c]{cc}%
h_{j} & 0\\
0 & h_{j}%
\end{array}
\right) \\
&&  =\frac{h_{j}}{h_{i}}\left(
\begin{array}
[c]{cc}%
\frac{df_{1}}{d\theta}(\gamma_{\beta}^{-1}(\tilde{\eta}_{j}^{-1}(p)) & 0\\
\frac{1}{\beta}\frac{df_{2}}{d\theta}(\gamma_{\beta}^{-1}(\tilde{\eta}%
_{j}^{-1}(p)) & \frac{df_{2}}{dx}(\gamma_{\beta}^{-1}(\tilde{\eta}_{j}%
^{-1}(p))
\end{array}
\right)  ,
\end{eqnarray*}
which in turn is arbitrarily close to $\frac{h_{j}}{h_{i}}diag(\frac{df_{1}%
}{d\theta},\frac{df_{2}}{dx}).$ This means that by using the artificial
rescaling $\gamma_{\beta}$ (without the actual need to apply it in practice
for our computer assisted proof), we can divide the region $\mathcal{V}$ into
a finite number of sets $U_{1},\ldots,U_{N}$ ($\mathcal{V}\subset\bigcup
_{i=1}^{N}U_{i}$), and verify cone conditions using interval matrices
$diag(\left[  \frac{df_{1}}{d\theta}(U_{i})\right]  ,\left[  \frac{df_{2}}%
{dx}(U_{i})\right]  )$ and applying Lemma \ref{lem:matrix-bounds-cc}. For our
proof we take $\mathbf{\gamma}_{0}=(a,b)=\left(  1,-1\right)  ,$ which means
that the quadratic form for our cones is simply%
\[
Q_{\mathbf{\gamma}_{0}}(\theta,x)=x^{2}-\theta^{2}.
\]
If we take $\mathbf{\gamma}_{1}=\left(  (1-\varepsilon),1\right)  $ for any
small parameter $\varepsilon>0$ then by choosing sufficiently large $\beta$
Assumption \ref{as:cones-setup} is satisfied (since any switch to new
coordinates is arbitrarily close to identity). This means that we can take
$\mathbf{\gamma}_{1}=\mathbf{\gamma}_{0}$, provided that all the inequalities
in our verification of cone conditions in the computer assisted proof are strict.

\subsubsection{Tools used for the proof\label{sec:tech-notes}}

Our proof has been conducted with the use of the CAPD library
(http://capd.ii.uj.edu.pl) developed by the Computer Assisted Proofs
in Dynamics group. We have used the multi-precision version of the library
running at 128 mantisa bits accuracy (which is approximately equivelent to
tracking 40 digits). The $C^{10}$ computations have been performed with
assistance of the Flexible Automatic Differentiation Package FADBAD++
(www.fadbad.com). The proof takes 16 seconds running on a 2.53 GHz laptop with
4GB of RAM.

\section{Final comments}

In this paper we have presented a version of a normally hyperbolic invariant manifold theorem,
which can be applied for rigorous-computer-assisted proofs. We have successfully applied our 
method to an example in which standard {\tt double} precision simulations brake down and produce 
false results. This demonstrates the strength of our method, that it can handle numerically difficult cases. It needs
to be noted that to apply our method we have used multiple precision for our computer assisted computations. 
For our proof we also needed to apply a high order method which relied on $C^{10}$ computations. 
We believe that it should be possible to devise a similar in spirit method, which would give proofs without multiple precision 
and using $C^1$ computations only. This will be the subject of our future work.

\section{Acknowledgements}

We would like to thank Tomasz Kapela for his assistence and comments regarding
the implementation of multi-precision in CAPD library. Our special thanks goes
to Daniel Wilczak for his suggestions, frequent discussions and for his
assistence with implementation of higher order computations in the CAPD library. The research of MC has been supported by the Polish State Ministry of Science and Information Technology grant N201 543238. 
The research of CS has been supported by grants MTM2006-05849/Consolider (Spain),
and CIRIT 2008SGR--67 (Catalonia).